\documentclass[a4paper,10pt]{amsart}
\usepackage[USenglish]{babel}
\usepackage[utf8]{inputenc}
\pdfoutput=1 

\usepackage[style=american]{csquotes}
\usepackage{amsmath,amsthm,amscd,amssymb}
\usepackage{MnSymbol}
\usepackage{smartref}
\usepackage{tikz}

\definecolor{color_cite}{RGB}{0,0,200}
\definecolor{color_link}{RGB}{200,0,0}
\definecolor{light_gray}{gray}{0.7}

\usepackage[hyperindex,breaklinks]{hyperref}
\hypersetup{
  colorlinks   = true,          
  urlcolor     = color_cite,          
  linkcolor    = color_link,          
  citecolor   = color_cite             
}

\usetikzlibrary{shapes}

\title{The intersection ring of matroids}
\author{Simon Hampe}

\address{Technische Universität Berlin\\
Institut für Mathematik, Sekretariat MA 6-2\\
Straße des 17. Juni 136, 10623 Berlin}
\email{hampe@math.tu-berlin.de}

\newcommand{\tpn}[1]{\mathbb{R}^{#1}/\textbf{1}}

\newcommand{\cone}{\textnormal{cone}}

\newcommand{\rank}{\textnormal{rank}}
\newcommand{\corank}{\textnormal{corank}}

\newcommand{\M}{\mathbb{M}}
\newcommand{\cyc}{\textnormal{cyc}}
\newcommand{\free}{\textnormal{free}}

\newcommand{\frank}{\textnormal{frk}}
\newcommand{\nty}{\textnormal{null}}
\newcommand{\lex}{\textnormal{lex}}
\newcommand{\cl}{\textnormal{cl}}
\newcommand{\TCP}{{\curly{C}_{\curly{Z}}}}
\newcommand{\nest}[1]{_{\left(#1\right)}}
\definecolor{light-gray}{gray}{0.8}

\renewcommand{\epsilon}{\varepsilon}

\tikzset{
  >=latex  
}

\theoremstyle{definition}
\newtheorem{defn}{Definition}[section]
\newtheorem{ex}[defn]{Example}

\newtheorem{remark}[defn]{Remark}

\newtheorem*{acknowledgement}{Acknowledgement}

\theoremstyle{plain}
\newtheorem{theorem}[defn]{Theorem}
\newtheorem{corollary}[defn]{Corollary}
\newtheorem{lemma}[defn]{Lemma}
\newtheorem{prop}[defn]{Proposition}

\newtheorem*{theoremnonr}{Theorem}

\newcommand{\abs}[1]{\left\lvert #1 \right\rvert} 
\newcommand{\gnrt}[1]{\left\langle #1 \right\rangle} 

\newcommand{\Z}{\mathbb{Z}}
\newcommand{\Q}{\mathbb{Q}}
\newcommand{\N}{\mathbb{N}}
\newcommand{\R}{\mathbb{R}}

\newcommand{\Star}{\textnormal{Star}}

\newcommand{\conv}{\textnormal{conv}}

\newcommand{\Ch}{\textnormal{Ch}}
\newcommand{\curly}[1]{\mathcal{#1}}

\newcommand{\wo}{\backslash} 

\begin{document}

\begin{abstract}
 We study a particular graded ring structure on the set of all loopfree matroids on a fixed labeled ground set, which occurs naturally in tropical geometry. The product is given by matroid intersection and the additive structure is defined by assigning to each matroid the indicator vector of its chains of flats. We show that this ring is generated in corank one, more precisely that any matroid can be written as a linear combination of products of corank one matroids. Moreover, we prove that a basis for the graded part of rank $r$ matroids is given by the set of nested matroids and that the total number of these is a Eulerian number. Derksen's $\mathcal{G}$-invariant then defines a $\Z$-linear map on this ring, which implies for example that the Tutte polynomial is linear on it as well. Finally we show that the ring is the cohomology ring of the toric variety of the permutohedron and thus fulfills Poincar\'e duality.
 \end{abstract}

\maketitle

\section{Introduction}\label{section_introduction}

In matroid theory, one can associate various algebraic invariants to individual matroids, such as, for example, the Tutte polynomial. However, one can also try to define algebraic structures on \emph{sets of matroids}. For example, Crapo and Schmitt study coalgebra and Hopf algebra structures on matroids \cite{sincidencehopf, csprimitivematroidminor, csuniquefactorization}. In \cite{abgwhomology}, the authors define homology groups of matroids. Recently, Giansiracusa and Giansiracusa \cite{gggrassmannalgebra} defined the Grassmann algebra of valuated matroids as an idempotent analogue of the classical Grassmann algebra. 

Our approach is inspired by intersection theory in tropical geometry. Tropical geometry can be seen as a combinatorial or polyhedral version of algebraic geometry, though it has ramifications into (among others) optimization, number theory, biological statistics and economics. We recommend the book by Maclagan and Sturmfels \cite{MaclaganSturmfelsBook} as an introduction to the subject.

Matroids have played an important role in tropical geometry ever since Sturmfels discovered that the tropicalization of a linear space only depends on the underlying matroid \cite{ssolving}. Speyer \cite{stropicallinear} generalized the concept of a tropical linear space to valuated matroids \cite{DressWenzel}. He also showed that there exists a product structure on valuated matroids by proving that the intersection product of two tropical linear spaces is a tropical linear space. It is this product structure that we wish to exploit, though we will only consider ordinary matroids, i.e.\ trivial valuations. Many geometric operations on tropical linear spaces have matroid-theoretic counterparts, which makes them very well-behaved and well-understood. Thus matroid theory has provided a very useful tool for tropical geometers (see for example \cite{smatroidintersection, frdiagonalintersection}).

Here we intend to go the opposite way: Use the geometric intuition behind tropical geometry to learn something about matroids. Indeed, the operations we will define to obtain the intersection ring of matroids might seem odd and unnatural to the matroid theorist. However, from the point of view of tropical geometry it is immediately clear that this ring is the right object to consider. It is the intersection ring of tropical linear spaces (with trivial valuation), where the operations are the obvious ones. As we will show, this ring exhibits a very rich structure that is tightly connected to the underlying matroids. In particular, it has the striking property that the $\mathcal{G}$-invariant induces a $\Z$-linear map on it. This invariant was introduced by Derksen in \cite{dsymmetricpolymatroids} as a generalization of various matroid invariants. In particular, all invariants which can be derived linearly from the $\curly{G}$-invariant automatically induce linear maps on the matroid intersection ring. This includes: The Tutte polynomial, all Tutte-Grothendieck invariants, the number of flats of fixed rank and the number of cyclic flats of fixed rank.
Derksen and Fink showed in \cite{dfvaluative} that $\curly{G}$ is the universal invariant for valuations on matroid polytopes. We will see in Section \ref{subsection_poincare_duality} that as a $\Z$-module the intersection ring of matroids is a quotient of the matroid polytope module they construct. 


The ring is naturally graded by corank and each graded piece is a free $\Z$-module whose dimension is a well-known combinatorial quantity: The Eulerian number $A_{r,n}$, which by definition counts permutations on $n$ elements with $r$ ascents. We show that it is also the number of \emph{nested matroids} of rank $r+1$ on $n$ labeled elements. To our knowledge, this is a new result -- only the number of \emph{isomorphism classes} of nested matroids had been determined so far.

The basic idea for defining the intersection ring of matroids is this (a formal definition will be made in Section \ref{section_preliminaries}): We identify each loopfree matroid of rank $r$ on $n$ labeled elements with its indicator vector of maximal chains of flats. We denote the $\Z$-module thus obtained by $\M_{r,n}$. We then define a product on $\M_n = \bigoplus_{r \geq 1} \M_{r,n}$ via \emph{matroid intersection} (which is the dual operation of matroid union): For any two matroids $M,N$, we set $M \cdot N := M \wedge N$ if the latter is loopfree, and $0$ otherwise. The product is extended to arbitrary linear combinations by distributivity. Phrased like this it is not even clear that this is well-defined. However, tropically this is clear: A linear combination of matroids corresponds to a linear combination of tropical linear spaces, which is an actual geometric object, to which an intersection product can be applied. Well-definedness then follows from the distributivity of the intersection product.

We now summarize our main results:

\begin{theoremnonr} \par\indent
\begin{itemize}
 \item With the operations defined above, $(\M_n,+,\cdot)$ is a commutative $\Z$-algebra with $1 = U_{n,n}$. It is graded by corank, i.e. 
 $$\M_{r,n} \cdot \M_{s,n} \subseteq \M_{r+s-n,n}$$
 (where $\M_{k,n} = 0$ if $k \leq 0$).
 Furthermore, it is generated in corank one: every matroid can be written as a linear combination of products of corank one matroids.
 \item  A basis of the free $\Z$-module $\M_{r,n}$ is given by the set of all loopfree nested matroids of rank $r$. The number of these matroids is the Eulerian number $A_{r-1,n}$.
 \item The $\mathcal{G}$-invariant induces a $\Z$-module homomorphism $\M_n \to \Z$ via
 $$M \mapsto \curly{G}(M) \textnormal{ for all matroids } M\enspace.$$
 \item $\M_n \cong A^*(X(\textnormal{Perm}_n))$, the cohomology ring of the toric variety corresponding to the normal fan of the permutohedron of order $n$. In particular, it fulfills Poincar\'e duality. That is, for every $1 \leq r \leq n$, the intersection product induces a perfect pairing
 $$\M_{r,n} \times \M_{n-r+1,n} \to \M_{1,n} \cong \Z\enspace.$$
 \item $\M_n$ is isomorphic to the subalgebra of McMullen's polytope algebra generated by all matroid polytopes.
\end{itemize}
\end{theoremnonr}

As the theorem suggests, a central role in our arguments is played by \emph{nested matroids}. These matroids have occurred under a variety of names, such as \emph{generalized Catalan matroids} \cite{bmlatticepathmatroids} or \emph{shifted matroids} \cite{acatalanmatroid} and seem to have been first defined by Crapo \cite{csingleelement}.
They are transversal matroids whose transversal presentation is a chain of sets. They are minor- and dual-closed, well-quasi-ordered and have an infinite list of excluded minors \cite{oprmatroidsdomains,bmlatticeofcyclic}. Their relevance to our approach stems from the fact that they can equivalently be characterized by the fact that their lattice of \emph{cyclic flats} is a chain. Cyclic flats are flats which are unions of circuits. They encode the full lattice of flats and nested matroids can thus be seen as the basic building blocks for constructing more complicated matroids.

The last part of the theorem above emphasizes the relevance of the multiplicative structure. The Eulerian numbers are symmetric, thus implying that $\M_{r,n}$ is isomorphic to $\M_{n-r+1,n}$ as a free $\Z$-module. This is surprising at first, as the isomorphism clearly cannot be induced by any obvious matroid operation such as taking duals (in fact, dualizing is not even a well-defined map on $\M_n$, as it may produce loops). However, the corresponding tropical linear spaces do have complementary dimension and so immediately suggest Poincar\'e duality.

The layout of the paper is as thus: In Section \ref{section_preliminaries}, we mainly introduce the central notions from matroid theory we will need -- in particular the notion of cyclic flats. A rigorous definition of the matroid intersection ring will be made. We also give a very brief summary of the most basic definitions from tropical geometry. In Section \ref{section_linear_combinations} we will only be concerned with the additive structure on $\M_{r,n}$, the product will not yet play a role. We show that the nested matroids are linearly independent and that they form, in fact, a basis. In Theorem \ref {theorem_basis_presentation} we will give an explicit formula for the representation of an arbitrary matroid in terms of this basis, which is based on the lattice of chains of cyclic flats of the matroid. To this end we introduce the notion of cyclic reductions, which are special weak maps. In \ref{subsection_tutte} we show that the $\mathcal{G}$-invariant and various other matroid invariants define $\Z$-linear maps on $\M_n$. Our proofs build on results of Bonin and Kung \cite{bkcatenary}. Section \ref{section_counting_nested} is quite short: We prove that the number of loopfree nested matroids of rank $r$ on $n$ labeled elements is the Eulerian number $A_{r-1,n}$. Section \ref{section_intersection_product} is dedicated to studying the multiplication on $\M_n$. In \ref{subsection_chain_matroids} we introduce the notion of chain products, which are special products of corank one matroids and we show that these are the same as nested matroids. In \ref{subsection_vanishing} we study in which cases the product of a matroid and a nested matroid vanishes. In the last part \ref{subsection_poincare_duality} we show that $\M_n$ is isomorphic to two familiar algebraic objects: The cohomology ring of a smooth, complete toric variety and the subring of the polytope algebra generated by matroid polytopes. The first result then immediately implies Poincar\'e duality (in fact, this is also a special case of a more general result by Adiprasito, Huh and Katz \cite{ahkhodgetheory}. Note also that Poincar\'e duality was known to hold for the full polytope algebra \cite{bstructurepolytopealgebra}). The second result shows that $\M_n$ is a quotient of the matroid polytope module considered by Derksen and Fink. The last part, Section \ref{section_outlook}, contains some suggestions for further research.


\begin{acknowledgement}
 The author was supported by DFG grant JO366/3-2, which is part of the DFG priority project SPP 1489 (\url{www.computeralgebra.de}). I would like to thank Michael Joswig and Benjamin Schröter for many helpful discussions and the anonymous referees for their very constructive comments.
 
 Many of the results in this paper were first discovered by computational means, using the author's software \texttt{a-tint} \cite{hatint}, which is an extension for \texttt{polymake} \cite{gjpolymake}. The latter is a software package for polyhedral and combinatorial computations. Functionality for computing in the intersection ring of matroids will be included in its next release.
\end{acknowledgement}

\section{Preliminaries}\label{section_preliminaries}

In this section, we mostly collect the definitions and results from matroid theory and combinatorics that are relevant to this paper. Part \ref{subsection_tropical} contains definitions from tropical geometry and illuminates the origin of the ring structure that we study here. We will assume familiarity with basic notions of matroid theory, for which we recommend \cite{omatroidtheory, wtheoryofmatroids} as references. 

\begin{defn} 
 All matroids are matroids on a labeled ground set, i.e.\ we do \emph{not} consider isomorphism classes of matroids. If not explicitly stated otherwise, all matroids are assumed to be loopfree.

 The complement of a set $A$ is written as $A^c$. We write $U_{r,E}$ for the uniform matroid of rank $r$ on the ground set $E$ and $U_{r,n}$ if $E = \{1,\dots,n\}$. 
 
 Let $M$ be a matroid on the ground set $E$. A \emph{flat} of $M$ is a set $F$ such that for any $x \notin F$, $\rank_M(F \cup x) = \rank_M(F) +1$. The set of flats of $M$ is denoted by $\curly{F}(M)$. The set of flats of rank $s$ is $\curly{F}_s(M)$. The \emph{closure} of a set $A$ in $M$, denoted by $\cl_M(A)$, is the smallest flat containing $A$. A \emph{spanning set} of $M$ is a set $A$ such that $\cl_M(A) = E$. The set of all spanning sets is $\curly{S}(M)$. The set of bases of $M$ is written $\curly{B}(M)$.
 
 The \emph{corank} of a matroid on $n$ elements is $\corank(M) = n - \rank(M)$. Similarly, the corank of a set $A$ in $M$ is $\corank_M(A) = \rank(M) - \rank_M(A)$. The \emph{nullity} of a set $A$ is $\nty_M(A) = \abs{A} - \rank_M(A)$. 
 
 The dual of $M$ is denoted by $M^*$. A matroid $N$ on $E$ is a \emph{quotient} of a matroid $M$ on $E$ if every flat of $N$ is also a flat of $M$. Equivalently, $N$ is a quotient of $M$ if and only if every circuit of $M$ is a union of circuits of $N$.
\end{defn}

\subsection{Posets, lattices and the Möbius function}\label{subsection_lattices}

\begin{defn}
Let $(\curly{P}, \leq)$ be a poset.
\begin{itemize}
 \item Let $x,y,z \in \curly{P}$. We say that $z$ is the \emph{join} of $x$ and $y$ in $\curly{P}$, written $z = x \vee_{\curly{P}} y$, if it is the unique minimal element that is larger than or equal to both $x$ and $y$. Similarly, $z$ is the \emph{meet} $x \wedge_{\curly{P}} y$ if it is the unique maximal element which is smaller than or equal to both $x$ and $y$.
 \item We call $\curly{P}$ a \emph{lattice} if for any two elements $x,y \in \curly{P}$ both join and meet exist.
\item We say that $\curly{P}$ is \emph{join-contractible} if there exists an $a \in \curly{P}$, such that $x \vee_{\curly{P}} a$ exists for all $x \in \curly{P}$.
\item For a poset $\curly{P}$ we denote by $\Ch(\curly{P})$ the poset of nonempty chains of elements of $\curly{P}$, partially ordered by inclusion. Here a chain means a set of pairwise comparable elements.
\item The \emph{order complex} $\Delta(\curly{P})$ of $\curly{P}$ is the simplicial complex on the vertex set $\curly{P}$ whose faces are all chains in $\curly{P}$.
\item We denote by $\hat{\curly{P}} := \curly{P} \cup \{\hat{0},\hat{1}\}$ the poset obtained from $\curly{P}$ by adjoining artificial minimal and maximal elements, so that $\hat{0} < x < \hat{1}$ for all $x \in \curly{P}$.
\item We say that an element $x$ of $\curly{P}$ \emph{covers} another element $y$ if $x > y$ and there is no $z$ such that $x > z > y$.
\end{itemize}
\end{defn}

\begin{ex}
 The basic example of a lattice in the context of matroids is the set of flats $\curly{F}(M)$ of a matroid $M$. A partial ordering is induced by set inclusion, the meet of two flats is their intersection and the join is the closure of the union.
\end{ex}

\begin{defn}
 Let $\curly{P}$ be a poset. We define the \emph{Möbius function} $\mu_{\curly{P}} : \curly{P} \times \curly{P} \to \Z$ in the following, recursive manner:
 \begin{itemize}
  \item $\mu_{\curly{P}}(x,y) = 0$ if $x \nleq y$.
  \item $\mu_{\curly{P}}(x,x) = 1$.
  \item Let $x \leq y$ and assume $\mu_{\curly{P}}(x,z)$ has been defined for all $x \leq z < y$. Then $$\mu_{\curly{P}}(x,y) = - \sum_{x \leq z < y} \mu_{\curly{P}}(x,z)\enspace.$$
 \end{itemize}
 We define the \emph{Möbius number} of $\curly{P}$ to be $\mu(\curly{P}) := \mu_{\hat{\curly{P}}}(\hat{0},\hat{1}).$
\end{defn}

There is an abundance of literature on the Möbius function and its various properties, see for example \cite{rfoundationscombinatorial, acombinatorialtheory} and the unpublished notes by Godsil \cite{gintroductionmoebius}. Our notation follows the latter. We will need the following well-known results (see also \cite{btopologicalmethods} for a more topological formulation):

\begin{lemma}[e.g.\ {\cite[Prop.\ 4.6]{acombinatorialtheory}}]\label{lemma_moebius}
 Let $\curly{P}$ be a poset and $x,z \in \curly{P}$. Then 
 $$\sum_{x \leq y \leq z} \mu_{\curly{P}}(y,z) = 0\enspace.$$
\end{lemma}

\begin{lemma}[e.g.\ {\cite[Cor.\ 10.13]{btopologicalmethods}}]\label{lemma_join_contractible}
 Let $\curly{P}$ be a join-contractible poset. Then
 $$\mu(\curly{P}) =0\enspace.$$
\end{lemma}

The next result is a consequence of the fact that the Möbius number of $\curly{P}$ is the reduced Euler characteristic of $\Delta(\curly{P})$ (e.g.\ \cite[Prop.\ 3.6]{rfoundationscombinatorial}) and that $\Delta(\curly{P})$ is homeo\-morphic to $\Delta(\textnormal{Ch}(\curly{P}))$ (the latter is the barycentric subdivision of the first).

\begin{prop}\label{prop_moebius_chains}
 Let $\curly{P}$ be a poset. Then
 $$\mu(\curly{P}) = \mu(\Ch(\curly{P}))\enspace.$$
\end{prop}

\subsection{Cyclic flats, transversal and nested matroids}\label{subsection_cyclic_flats}

\begin{defn}\label{def_cyclic_flats}
 For each flat $F$ of a matroid $M$, we define the \emph{cyclic part} of $F$ to be 
 $$\cyc_M(F) := \bigcup_{\substack{C \textnormal{ circuit of }M\\ C \subseteq F}} C\enspace.$$
 This is again a flat of $M$. We write $\free_M(F) := F \wo \cyc_M(F)$ and $\frank_M(F) = \abs{\free(F)}$. We note some obvious properties:
 \begin{itemize}
  \item $M_{\mid F} = M_{\mid \cyc(F)} \oplus U_{\frank(F),\free(F)}$ and $\rank_M(F) = \rank_M(\cyc_M(F)) + \frank(F)$.
  \item Every set $\cyc_M(F) \subseteq A \subseteq F$ is again a flat of $M$.
 \end{itemize}
 We call a flat $F$ cyclic if $F = \cyc_M(F)$. The set of all cyclic flats of $M$ is denoted by $\curly{Z}(M)$. 
 
 We say that $M$ is \emph{nested} if $\curly{Z}(M)$ is a chain of sets.
\end{defn}

\begin{remark}\label{remark_cyclic_construction}
 Thomas Brylawski pointed out in \cite{baffinerepresentation} that knowing the ground set, all cyclic flats and their ranks is sufficient to determine the whole matroid. The actual construction will be relevant to some of our arguments, so we recall it here.
 
 We partition the set of all flats $\curly{F}(M)$ into the sets
 $$\curly{F}(M)_{s,m} := \{F \in \curly{F}(M); \rank(F) = s \textnormal{ and } \frank(F) = m\}\enspace.$$
 We then recursively build up $\curly{F}(M)$ in the following manner:
  \begin{enumerate}
  \item $\curly{F}(M)_{s,0}$ is the set of all cyclic flats of rank $s$.
  \item For $m > 0$, $\curly{F}(M)_{s,m}$ is the set of all $F \cup \{p\}$, such that:
  \begin{itemize}
   \item $F \in \curly{F}(M)_{s-1,m-1}$ and $p \notin F$.
   \item $F \cup \{p\}$ is not contained in any $G$, where $G \in \curly{F}(M)_{s,m'}$ and $m' < m$.
  \end{itemize}
 \end{enumerate}
 \end{remark}
 
 Bonin and de Mier later proved that this information in fact provides a cryptomorphic characterization of matroids. 
 
 \begin{theorem}[{\cite[Theorem 3.2]{bmlatticeofcyclic}}]\label{thm_cyclic_axioms}
  Let $\curly{Z}$ be a collection of subsets of $E$ and $r$ an integer-valued function on $\curly{Z}$. Then $\curly{Z}$ is the collection of cyclic flats of a matroid $M$ and $r$ the restriction of the rank function on $M$ if and only if the following hold:
\begin{itemize}
 \item[(Z0)] $\curly{Z}$ is a lattice under inclusion.
 \item[(Z1)] $r(0_{\curly{Z}}) = 0$, where $0_{\curly{Z}}$ is the minimal element of $\curly{Z}$.
 \item[(Z2)] $0 < r(Y) - r(X) < \abs{Y - X}$ for all sets $X,Y$ in $\curly{Z}$ with $X \subsetneq Y$.
 \item[(Z3)] For all sets $X,Y$ in $\curly{Z}$,
  $$r(X) + r(Y) \geq r(X \vee_{\curly{Z}} Y) + r(X \wedge_{\curly{Z}} Y) + \abs{ (X \cap Y) - (X \wedge_{\curly{Z}} Y)}\enspace.$$
\end{itemize}
\end{theorem}

\begin{remark}
Note that in a matroid, $X \vee_{\curly{Z}(M)} Y = \cl_M(X \cup Y)$ is the usual join of flats, but $X \wedge_{\curly{Z}(M)} Y$ is the union of all circuits contained in $X \cap Y$ and can in general be strictly smaller than the flat $X \cap Y$. Also, $0_{\curly{Z}(M)}$ is the set of loops and $1_{\curly{Z}(M)}$ is the union of all circuits and thus the complement of the coloops of $M$.
\end{remark}

\begin{ex}\label{ex_cyclic_flats}\par\indent
\begin{enumerate}
 \item We define a rank two matroid $M$ on $E = \{1,\dots,4\}$ via its lattice of flats:
 $$\curly{F}(M) := \{ \emptyset, \{1,4\},\{2,3\}, E\}\enspace.$$
 Its circuits are $\{1,4\}$ and $\{2,3\}$, so $\curly{Z}(M) = \curly{F}(M).$ In particular, $M$ is not nested. 
 \item The circuits of the uniform matroid $U_{r,n}$ are all the sets of size $r+1$. Its flats are
 $\curly{F}(U_{r,n}) = \{ F \subseteq E,\; \abs{F} < r\} \cup \{E\}.$
 Hence $\curly{Z}(U_{r,n}) = \{\emptyset, E\}$, so any uniform matroid is nested.
\end{enumerate}
\end{ex}

\begin{defn}
 Let $\curly{A} := (A_1,\dots,A_m)$ be subsets of $E$ (which need not be distinct). A \emph{partial transversal} of $E$ is a subset $S \subseteq E$ such that there is a bijection $\psi: J \to S$ from a set $J \subseteq [m]$ fulfilling $\psi(j) \in A_j$ for all $j \in J$. The set of all partial transversals forms the set of independent sets of a matroid (see for example \cite[Theorem 1.6.2]{omatroidtheory}), which we denote by $M[A_1,\dots,A_m]$. We call a matroid of this form a \emph{transversal} matroid.
\end{defn}

\begin{remark}\label{remark_transversal_presentations}
 The \emph{transversal presentation}, i.e.\ the choice of set system for a trans\-versal matroid is not unique. However, the following holds \cite{bwsomeresults}:
 \begin{itemize}
  \item Every rank $s$ transversal matroid has a presentation with $s$ sets $A_1,\dots,A_s$. More precisely, if $M = M[A_1,\dots,A_k]$ has a basis $B$ which is a transversal of $A_{i_1},\dots,A_{i_s}$, then $M = M[A_{i_1},\dots,A_{i_s}]$.
  \item Every rank $s$ transversal matroid $M$ has a unique \emph{maximal presentation}, i.e.\ one cannot increase any of the sets $A_i$ without changing the matroid. This maximal presentation is constructed as follows:
  \begin{enumerate}
   \item Assume $M = M[A_1',\dots,A_s']$ has a presentation with $s$ sets.
   \item For each $i = 1,\dots,s$, replace $A_i'$ by
   $$A_i := A_i' \cup R_i\enspace,$$
   where $R_i$ is the set of coloops of $M_{\mid A_i'^c}$.
  \end{enumerate}
 In particular, a presentation $M[A_1,\dots,A_s]$ is maximal if and only if $M_{\mid A_i^c}$ has no coloops.
 \item The restriction of a transversal matroid is again transversal. To be precise, if $T \subseteq E$, we have
 $$M[A_1,\dots,A_s]_{\mid T} = M[A_1 \cap T,\dots,A_s \cap T]\enspace.$$
 \end{itemize}
\end{remark}

The survey \cite{bintroductiontransversal} is a good introduction to the subject of transversal matroids. They are relevant to our discussion mainly because of the following result, a proof of which can for example be found in \cite{oprmatroidsdomains}:

\begin{theorem}\label{theorem_nested_matroids}
 A matroid $M$ is nested if and only if $M = M[A_1,\dots,A_k]$ for a chain of sets $A_1 \subseteq \dots \subseteq A_k$.
\end{theorem}

\begin{ex}\label{ex_transversal}
 The matroid $M$ from Example \ref{ex_cyclic_flats},(1) has the maximal presentation $M = M[\{1,4\}, \{2,3\}]$. The uniform matroid $U_{r,n}$ is the transversal matroid $M[E,\dots,E]$, where $E = [n]$ occurs $r$ times.
\end{ex}

\subsection{The intersection ring of matroids}\label{subsection_intersection_ring}

The notion of matroid intersection seems to have gotten very little attention from matroid theorists. It is the dual operation of a much more actively studied object, the matroid union (see for example \cite[Chapter 11.3]{omatroidtheory} or \cite[Chapter 7.6]{wtheoryofmatroids}):

\begin{defn}
 Let $M,N$ be matroids on a ground set $E$. The \emph{union} $M \vee N$ is the matroid on $E$, whose independent sets are
 $$\curly{I}_{M \vee N} := \{ I \cup J; I \in \curly{I}_M, J \in \curly{I}_N\}\enspace.$$
 The \emph{intersection} of $M$ and $N$ is then defined as
 $$M \wedge N = (M^* \vee N^*)^*\enspace.$$
\end{defn}

\begin{remark}
We note a few properties of matroid intersection:

\begin{enumerate}
 \item It is known that both $M$ and $N$ are quotients of $M \vee N$. By duality it follows that $M \wedge N$ is a quotient of both $M$ and $N$ (see \cite[Chapter 7.3]{omatroidtheory}).
 \item The spanning sets of $M \wedge N$ are given by:
 $$\curly{S}(M \wedge N) = \{ S \cap T; S \in \curly{S}(M), T \in \curly{S}(N)\}\enspace.$$
 \item As bases are the minimal spanning sets of a matroid and since $M \wedge N$ is a quotient of both $M$ and $N$, this implies
 $$\curly{B}(M \wedge N) = \{B \cap B'; B \in \curly{B}(M), B' \in \curly{B}(N), \abs{B \cap B'} = \rank(M \wedge N)\}\enspace.$$
 \item Matroid intersection commutes with contraction, i.e. $$(M \wedge N)/A = (M/A) \wedge (N/A)\enspace.$$
\end{enumerate}
\end{remark}

\begin{ex}\label{ex_intersection}
 Let $M$ be any matroid on $[n]$ of rank $r > 1$. It is an easy combinatorial exercise to see that the bases of $M \wedge U_{n-1,n}$ are all the sets of the form $\{B' \subseteq B; B \textnormal{ a basis of }M\textnormal{ and }\abs{B'} = r - 1\}$. This is also called the \emph{truncation} $T(M)$ of $M$. Inductively, we see that $U_{n-k,n}$ is the $k$-fold intersection of $U_{n-1,n}$ with itself for any $k$, so $M \wedge U_{n-k,n}$ is just the $k$-fold truncation of $M$. In particular, $M \wedge U_{n-r+1,n} = U_{1,1}$ for all loopfree matroids $M$. 
\end{ex}

\begin{defn}\label{defn_matroid_ring}
 For $1 \leq r \leq n = \abs{E}$, let $\mathfrak{C}_{r,n}$ be the set of all chains of sets of length $r$, i.e.\ which are of the form
 $$\emptyset \subsetneq F_1 \subsetneq \dots \subsetneq F_r = E\enspace.$$
 We denote by $V_{r,n} = \Z^{\mathfrak{C}_{r,n}}$ the free $\Z$-module whose coordinates are indexed by the elements of $\mathfrak{C}_{r,n}$. 
 
 Let $\M_{r,n}^{\free}$ be the free $\Z$-module with generators the set of all loopfree matroids of rank $r$ on the ground set $\{1,\dots,n\}$. We define a homomorphism 
 $$\Phi_{r,n}: \M_{r,n}^{\free} \to V_{r,n}; M \mapsto v_M\enspace,$$
 where for each chain $\curly{C}$ we define
 $$(v_M)_{\curly{C}} := \begin{cases}
                        1, &\textnormal{if } \curly{C} \textnormal{ is a chain of flats in }M\\
                        0, &\textnormal{otherwise.}
                       \end{cases}
$$
The \emph{intersection ring of matroids on} $[n]$ is the $\Z$-module
$$\M_n = \bigoplus_{r=1}^n \M_{r,n}\enspace,$$
with $\M_{r,n} = \M_{r,n}^{\free} / \ker \Phi_{r,n}.$ It becomes a ring with the product defined by
$$M \cdot N := \begin{cases}
                M \wedge N, &\textnormal{if } M \wedge N \textnormal{ is loopfree}\\
                0, &\textnormal{otherwise,}
               \end{cases}$$
extended to linear combinations of matroids via distributivity.
\end{defn}

\begin{remark}
 It is not at all obvious that this is well-defined, i.e.\ that the definition of the product is compatible with the additive structure on $\M_n$. However, it follows implicitly from the fact that it is a tropical intersection product (see Remark \ref{remark_speyer_isomorphism}).
\end{remark}

\begin{ex}\label{ex_modules}
 The modules $\M_{r,n}$ are easy to write down for $r \in \{1,n\}$. In both cases there is only one loopfree matroid of rank $r$ on $n$ elements: $U_{1,n}$ and $U_{n,n}$, respectively. Hence $\M_{1,n} \cong \M_{n,n} \cong \Z$. We will postpone a nontrivial calculation until Example \ref{ex_relation}, where we can make use of the geometric intuition of tropical cycles.
\end{ex}

\subsection{Tropical geometry}\label{subsection_tropical}

In this section we will only sketch the most important definitions from tropical geometry. For a more in-depth account we recommend the book by Maclagan and Sturmfels \cite{MaclaganSturmfelsBook} and the book in progress by Mikhalkin and Rau \cite{mrtropicalgeometry}. Note that we use the $\min$-convention in our definition of matroid fans and that all coordinates are tropical projective coordinates in $\tpn{n} := \R^n/\gnrt{(1,\dots,1)}$.

\begin{defn}
 A \emph{tropical cycle} $(X,\omega_X)$ is a pure-dimensional, rational polyhedral complex $X$ in $\tpn{n}$ together with a function $\omega_X: X^{\max} \to \Z$ defined on its maximal cells that fulfills a certain \emph{balancing condition}: For a cone $\sigma$, we write $V_\sigma = \gnrt{a-b; a,b \in \sigma}$ and $\Lambda_\sigma = V_\sigma \cap \Z^n/\textbf{1}$. Then at every codimension one face $\tau$ of $X$ we must have
 $$\sum_{\sigma > \tau} \omega_X(\sigma) u_{\sigma/\tau} = 0 \;(\textnormal{mod } V_\tau)\enspace,$$
 where $u_{\sigma/\tau}$ is the primitive generator of the group $\Lambda_\sigma / \Lambda_\tau \cong \Z$ pointing towards $\sigma$.
 
 The \emph{support} of $X$ is the set $\abs{X} := \bigcup_{\sigma \in X^{\max}: \omega_X(\sigma) \neq 0} \sigma$.
 
 We will consider two tropical cycles to be the same if their supports have a common refinement respecting both weight functions (in particular, cells of weight zero are considered irrelevant).
 
 For any point $p \in \abs{X}$, we define the \emph{Star} of $X$ at $p$ to be the fan 
 $$\Star_X(p) := \{ \R_{\geq 0} \cdot (\sigma - p); p \in \sigma\}\enspace,$$
 with weight function $\omega_\Star(\R_{\geq 0} (\sigma - p)) = \omega_X(\sigma)$. It is easy to see that this is a tropical cycle.
 
 
 The \emph{sum} of two $k$-dimensional cycles $X,Y$ in $\tpn{n}$ is defined as follows: Choose refinements of $X$ and $Y$ such that they are defined on the same polyhedral structure (possibly defining some weights to be 0). Then 
 $$X + Y := X \cup Y \textnormal{ with weight function }\omega_{X + Y} := \omega_X + \omega_Y\enspace.$$ This operation makes the set of all $k$-dimensional cycles in $\tpn{n}$ into a group, which we denote by $Z_k(\tpn{n})$.
\end{defn}

\begin{defn}
 Let $M$ be a loopfree matroid on $[n]$. We define its \emph{matroid fan} $B(M)$ to be the fan in $\tpn{n}$ consisting of all cones
 $$\cone(\curly{C}) = \left\{ \sum_{i=1}^k \alpha_i v_{F_i}; \alpha_i \geq 0\right\}\enspace,$$
 where $\curly{C} = (F_1 \subsetneq \dots \subsetneq F_k)$ is a chain of flats in $M$ and $v_F = \sum_{i \in F} e_i \in \tpn{n}$.
\end{defn}

\begin{remark}
 This is a polyhedral fan of pure dimension $\rank(M)-1$. If one equips all maximal cells with weight 1, it becomes a tropical cycle. The notation $B(M)$ is in honor of George Bergman \cite{blogarithmiclimit}, who studied objects like these as logarithmic limit sets of algebraic varieties. The polyhedral structure given above was discovered by Ardila and Klivans \cite{akbergman}. The interested reader can find more information on matroid fans in the context of tropical geometry in \cite[Chapter 4]{MaclaganSturmfelsBook}. Note that each matroid fan is a subfan of the fan whose set of maximal cones is $\{\cone(\curly{C})\}$, where $\curly{C}$ runs over \emph{all} chains $\emptyset \subsetneq F_1 \subsetneq \cdots \subsetneq F_r = E$. Thus we can also identify a matroid fan with its indicator vector of chains $v_M \in V_{r,n}$ and the sum of two matroid fan cycles is just the sum of the indicator vectors. Hence we have:
\end{remark}

\begin{prop}
 $\M_{r,n}$ is isomorphic to the subgroup of $Z_{r-1}(\tpn{n})$ generated by matroid fans.
\end{prop}

Note that under this identification, a linear combination of matroids now has actual geometric meaning: It is the tropical cycle sum of the corresponding matroid fans. This makes dealing with the product much easier in this context, since we will define it on arbitrary tropical cycles. 

\begin{ex}\label{ex_relation}
 We encounter the first nontrivial linear relations for $n = 4$ and $r = 2$. We define four matroids on $E = \{1,\dots,4\}$ in terms of their flats:
 \begin{align*}
  \curly{F}(M_1) &:= \{\emptyset, \{1\},\{2\},\{3\},\{4\},E\}\enspace,\\
  \curly{F}(M_2) &:= \{\emptyset, \{1,4\},\{2,3\},E\}\enspace,\\
  \curly{F}(M_3) &:= \{\emptyset, \{1,4\},\{2\},\{3\},E\}\enspace,\\
  \curly{F}(M_4) &:= \{\emptyset, \{1\},\{2,3\},\{4\},E\}\enspace.
 \end{align*}
 Then one sees easily that $M_1 + M_2 = M_3 + M_4$ in $\M_{2,4}$. The corresponding tropical cycles are depicted in Figure \ref{figure_tropical_sum}.
 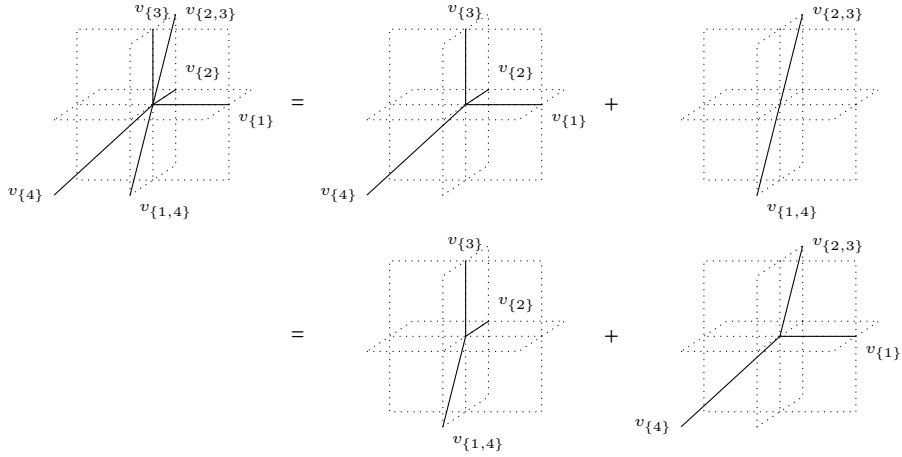
\begin{figure}[ht]
  \begin{tikzpicture}[x={(1cm,0cm)}, y={(0.3cm,0.2cm)},z={(0cm,1cm)}]
   \matrix[column sep = -1pt] {
    \draw (0,0,0) -- (1,0,0) node[below right]{\tiny $v_{\{1\}}$};
    \draw (0,0,0) -- (0,1,0) node [above right] {\tiny $v_{\{2\}}$};
    \draw (0,0,0) -- (0,0,1) node [above] {\tiny $v_{\{3\}}$};
    \draw (0,0,0) -- (-1,-1,-1) node [left] {\tiny $v_{\{4\}}$};
    \draw (0,0,0) -- (0,1,1) node[right] {\tiny $v_{\{2,3\}}$};
    \draw (0,0,0) -- (0,-1,-1) node[below right] {\tiny $v_{\{1,4\}}$}; 
    \draw[dotted] (1,0,0) -- (1,0,1) -- (-1,0,1) -- (-1,0,-1) -- (1,0,-1) -- (1,0,0);
    \draw[dotted] (-1,-1,0) -- (1,-1,0) -- (1,1,0) -- (-1,1,0) -- (-1,-1,0);
    \draw[dotted] (0,1,0) -- (0,1,1) -- (0,-1,1) -- (0,-1,-1) -- (0,1,-1) -- (0,1,0);
    \draw[dotted] (1,0,0) -- (-1,0,0);
    \draw[dotted] (0,1,0) -- (0,-1,0);
    \draw[dotted] (0,0,1) -- (0,0,-1);
    &
    \draw (0,0) node {$=$}; &
    \draw (0,0,0) -- (1,0,0) node[below right]{\tiny $v_{\{1\}}$};
    \draw (0,0,0) -- (0,1,0) node [above right] {\tiny $v_{\{2\}}$};
    \draw (0,0,0) -- (0,0,1) node [above] {\tiny $v_{\{3\}}$};
    \draw (0,0,0) -- (-1,-1,-1) node [left] {\tiny $v_{\{4\}}$}; 
    \draw[dotted] (1,0,0) -- (1,0,1) -- (-1,0,1) -- (-1,0,-1) -- (1,0,-1) -- (1,0,0);
    \draw[dotted] (-1,-1,0) -- (1,-1,0) -- (1,1,0) -- (-1,1,0) -- (-1,-1,0);
    \draw[dotted] (0,1,0) -- (0,1,1) -- (0,-1,1) -- (0,-1,-1) -- (0,1,-1) -- (0,1,0);
    \draw[dotted] (1,0,0) -- (-1,0,0);
    \draw[dotted] (0,1,0) -- (0,-1,0);
    \draw[dotted] (0,0,1) -- (0,0,-1);
    &
    \draw (0,0) node {$+$}; &
    \draw (0,0,0) -- (0,1,1) node[right] {\tiny $v_{\{2,3\}}$};
    \draw (0,0,0) -- (0,-1,-1) node[below right] {\tiny $v_{\{1,4\}}$};
    \draw[dotted] (1,0,0) -- (1,0,1) -- (-1,0,1) -- (-1,0,-1) -- (1,0,-1) -- (1,0,0);
    \draw[dotted] (-1,-1,0) -- (1,-1,0) -- (1,1,0) -- (-1,1,0) -- (-1,-1,0);
    \draw[dotted] (0,1,0) -- (0,1,1) -- (0,-1,1) -- (0,-1,-1) -- (0,1,-1) -- (0,1,0);
    \draw[dotted] (1,0,0) -- (-1,0,0);
    \draw[dotted] (0,1,0) -- (0,-1,0);
    \draw[dotted] (0,0,1) -- (0,0,-1);
    \\
    &
    \draw (0,0) node {$=$}; &
    \draw (0,0,0) -- (0,1,0) node [above right] {\tiny $v_{\{2\}}$};
    \draw (0,0,0) -- (0,0,1) node [above] {\tiny $v_{\{3\}}$};
    \draw (0,0,0) -- (0,-1,-1) node[below right] {\tiny $v_{\{1,4\}}$}; 
    \draw[dotted] (1,0,0) -- (1,0,1) -- (-1,0,1) -- (-1,0,-1) -- (1,0,-1) -- (1,0,0);
    \draw[dotted] (-1,-1,0) -- (1,-1,0) -- (1,1,0) -- (-1,1,0) -- (-1,-1,0);
    \draw[dotted] (0,1,0) -- (0,1,1) -- (0,-1,1) -- (0,-1,-1) -- (0,1,-1) -- (0,1,0);
    \draw[dotted] (1,0,0) -- (-1,0,0);
    \draw[dotted] (0,1,0) -- (0,-1,0);
    \draw[dotted] (0,0,1) -- (0,0,-1);
    &
    \draw (0,0) node {$+$}; &
    \draw (0,0,0) -- (1,0,0) node[below right]{\tiny $v_{\{1\}}$};
    \draw (0,0,0) -- (-1,-1,-1) node [left] {\tiny $v_{\{4\}}$};
    \draw (0,0,0) -- (0,1,1) node[right] {\tiny $v_{\{2,3\}}$};
    \draw[dotted] (1,0,0) -- (1,0,1) -- (-1,0,1) -- (-1,0,-1) -- (1,0,-1) -- (1,0,0);
    \draw[dotted] (-1,-1,0) -- (1,-1,0) -- (1,1,0) -- (-1,1,0) -- (-1,-1,0);
    \draw[dotted] (0,1,0) -- (0,1,1) -- (0,-1,1) -- (0,-1,-1) -- (0,1,-1) -- (0,1,0);
    \draw[dotted] (1,0,0) -- (-1,0,0);
    \draw[dotted] (0,1,0) -- (0,-1,0);
    \draw[dotted] (0,0,1) -- (0,0,-1);
    \\
   };
  \end{tikzpicture}
  \caption{A tropical cycle which can be written as the sum of matroid fans in different ways. Note that all weights are 1 and that we draw a vector in $\tpn{4} \cong \R^3$ by choosing the representative whose last coordinate is 0. E.g.\ $v_{\{4\}} = (-1,-1,-1)$. The dotted lines just indicate coordinate hyperplanes.}\label{figure_tropical_sum}
 \end{figure}

\end{ex}

%
%

%
%

 There is a notion of \emph{intersection product} $X \cdot Y$ of two tropical cycles $X,Y$, which makes $Z_{n-1} = \bigoplus_{k \in \Z} Z_k(\tpn{n})$ into a ring (we set $Z_k = 0$ for $k \notin \{0,\dots,n-1\}$). There are various equivalent definitions of this product (\cite{ mtropicalapplications, arfirststeps, jystableintersection}). For the sake of legibility, we adopt the following description from \cite{jystableintersection} of the set $X \cdot Y$ and omit the definition of the weights. Since we will only consider products of matroid fans, which -- as we shall shortly see -- are again matroid fans, all occurring weights are one anyway.
 \begin{defn}
  Let $X,Y$ be two tropical cycles in $\tpn{n}$. Then 
  $$\abs{X\cdot Y} = \{ p \in \abs{X} \cap \abs{Y}; \dim( \Star_X(p) \boxplus \Star_Y(p)) = n-1\}\enspace,$$
  where $A \boxplus B = \{a + b; a \in A, b \in B\}$ is the Minkowski sum of sets.
%
 \end{defn}
 
 \begin{ex} 
  Let $X = Y = B(U_{2,3})$. This is the one-dimensional fan in $\tpn{3}$, whose three rays are spanned by the vectors $e_i, i=1,\dots,3$. We wish to compute the support of $X \cdot Y$. At every point $p \neq 0$, $\Star_X(p) = \Star_Y(p) = \Star_X(p) \boxplus \Star_Y(p)$ is an actual line. For $p = 0$, we have $\Star_X(p) \boxplus \Star_Y(p) = X \boxplus Y = \tpn{3}$, so $\abs{X \cdot Y} = \{0\}$, which is in fact the support of $B(U_{1,3} = U_{2,3} \wedge U_{2,3})$ (see also Figure \ref{figure_intersection}).
  \begin{figure}[ht]
   \begin{tikzpicture}
    \draw (0,0) -- (2,0) node[right] {$v_{\{1\}}$};
    \draw (0,0) -- (0,2) node[right] {$v_{\{2\}}$};
    \draw (0,0) -- (-1.5,-1.5) node[left] {$v_{\{3}\}$};
    \fill[black] (0,1) circle (2pt) node[right=5pt] {
    \begin{tikzpicture}[scale=0.25]
    \pgfgettransform\mytrafo
    \matrix[execute at begin cell=\pgfsettransform\mytrafo]{
    \draw (0,-1) -- (0,1); 
    \fill[black] (0,0) circle (5pt); &
    \draw (0,0) node {$\boxplus$}; &
    \draw (0,-1) -- (0,1); 
    \fill[black] (0,0) circle (5pt); &
    \draw (0,0) node {$=$}; &
    \draw (0,-1) -- (0,1); 
    \\
    };
    \end{tikzpicture}
    };
    \fill[black] (0,0) circle (2pt);
    \draw (0,-.5) node[right=1pt] {
    \begin{tikzpicture}[scale=0.25]
    \pgfgettransform\mytrafo
    \matrix[execute at begin cell=\pgfsettransform\mytrafo]{
    \draw (0,0) -- (0,1); 
    \draw (0,0) -- (1,0);
    \draw (0,0) -- (-.8,-.8);
    \fill[black] (0,0) circle (5pt); &
    \draw (0,0) node {$\boxplus$}; &
    \draw (0,0) -- (0,1); 
    \draw (0,0) -- (1,0);
    \draw (0,0) -- (-.8,-.8); 
    \fill[black] (0,0) circle (5pt); &
    \draw (0,0) node {$=$}; &
    \fill[light_gray] (1,1) -- (1,-1) -- (-1,-1) -- (-1,1) -- (1,1);
    \\
    };
    \end{tikzpicture}
    };
   \end{tikzpicture}
   \caption{The self-intersection of this tropical cycle is just the origin. Again, we draw points in $\tpn{3} \cong \R^2$ by choosing the representative whose last coordinate is zero.}\label{figure_intersection}
  \end{figure}
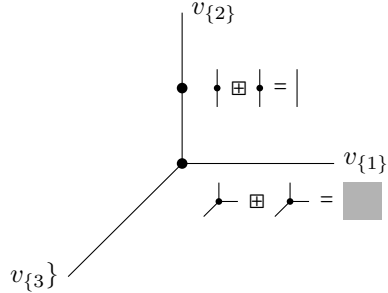
 \end{ex}

 \begin{remark}\label{remark_intersection_product}
$(Z_{n-1},+,\cdot)$, with cycle sum as addition and intersection product as multiplication is a $\Z$-algebra, graded by codimension (see \cite{arfirststeps} for details) and with multiplicative neutral element $1_{Z_{n-1}} = \tpn{n} = B(U_{n,n})$. In particular, we have 
  \begin{align*}
  Z_k \cdot Z_l &\subseteq Z_{k+l-n}\enspace,\\
  X \cdot (Y + Z) &= X \cdot Y + X \cdot Z\enspace.
  \end{align*}
\end{remark}

\begin{remark}\label{remark_speyer_isomorphism}
 David Speyer proved in \cite[Theorem 4.11]{stropicallinear} that $$B(M)\cdot B(N) = B(M \cdot N)$$ (where $B(0) = 0$, obviously). He proved this in the more general context of valuated matroids and tropical linear spaces. Together with remark \ref{remark_intersection_product}, this already implies the well-definedness of the product on $\M_n$ and that the resulting ring is a familiar object to tropical geometers:
\end{remark}

\begin{prop}
 $\M_n$ is isomorphic to the subring of the intersection ring of tropical cycles in $\tpn{n}$ which is generated by matroid fans. Hence it is a commutative ring with multiplicative neutral element $U_{n,n}$ and it is graded by corank.
\end{prop}

\section{Linear combinations of nested matroids}\label{section_linear_combinations}

This Section is dedicated to proving that the nested matroids of rank $r$ form a basis of $\M_{r,n}$. We will first prove linear independence in a more or less elementary manner. To show that each matroid can be written as a linear combination of nested matroids -- for which we will give an explicit formula -- requires more work.

\subsection{Linear independence}\label{subsection_linear_independence}

\begin{remark}\label{remark_cyclic_cardinality}
 By Remark \ref{remark_cyclic_construction} a nested matroid $M$ is uniquely determined by the list of tuples $$(Z_0 = \emptyset,r_0 =0),\dots,(Z_k,r_k),$$
 where $Z_0 \subsetneq \dots \subsetneq Z_k$ is the chain of cyclic flats of $M$ and $r_i = \rank_M(Z_i)$. We note some obvious properties of this data:
 \begin{align*}
  \nty_M(Z_k) &= \corank(M)\textnormal{, as $Z_k$ is the complement of the coloops of $M$}.\\
  \nty_M(Z_i) &< \nty_M(Z_j)\textnormal{ for all $i < j$, due to (Z2) of Theorem }\ref{thm_cyclic_axioms}.
 \end{align*}

\end{remark}

\begin{prop}
 The set of loopfree nested matroids of rank $r$ is linearly independent in $\M_{r,n}$.
 \begin{proof}
  For a nested matroid $M$ given by a tuple $$(Z_0 = \emptyset,r_0 := 0) , (Z_1,r_1), \dots, (Z_k, r_k)$$ as above, we define its \emph{gap measure} to be the tuple
  $$\gamma(M) := (d_i)_{i=1,\dots,r} \in \N^r, \textnormal{ where } d_i := \begin{cases}
                                                                            r_i - r_{i-1}, &\textnormal{if } i \leq k\\
									    0, &\textnormal{otherwise}.
                                                                           \end{cases}
$$

    Let $\curly{N}$ be the set of all loopfree nested matroids of rank $r$ on $[n]$. Assume there is a linear relation
    $$\sum_{M \in \curly{N}} a_M M = 0$$
    in $\M_{r,n}$. We will show by lexicographic induction on $\gamma(M)$ that $a_M =0$ for all $M$.
    
    If $\gamma(M) = (1,\dots,1)$, we have $r_i = i$ for all $i$. We complete the chain of cyclic flats to a chain of flats of $M$ of length $r$ in an arbitrary manner. We claim that $M$ is the only nested matroid containing this chain. In particular $a_M$ must be 0.
    
    Assume $N$ is a nested matroid of rank $r$ whose lattice of flats contains this chain. The empty set must be cyclic in $N$, as the matroid is loopfree. Inductively, we assume $Z_0,\dots,Z_{j-1}$ are cyclic in $N$. If $Z_j$ is noncyclic, it contains the cyclic flat $G := \cyc_N(Z_j)$ with $\abs{Z_j \wo G} = \rank_N(Z_j) - \rank_N(G)$. However, as $r_i = i$ for all $i$, we know that all the cyclic flats of smaller rank are $Z_0,\dots,Z_{j-1}$ and since $Z_j$ is cyclic in $M$ we know that $\abs{Z_j \wo Z_i} > \rank(Z_j) - \rank(Z_i)$ for any $i < j$. Hence $Z_j$ must be cyclic for all $j$. As $\nty_N(Z_k) = \nty_M(Z_k) = \corank(M) = \corank(N)$, $N$ can contain no cyclic flat larger than $Z_k$ and we conclude that $N = M$.
    
    If $\gamma(M) =: (d_1,\dots,d_r) > (0,\dots,0)$, we again complete the chain of cyclic flats to a chain of flats of $M$ of length $r$ in an arbitrary manner. We want to show that any nested matroid $N \neq M$ containing that chain must fulfill $\gamma(N) <_\lex \gamma(M)$, so by induction $a_N = 0$ and thus finally also $a_M = 0$.
    
    So let $N$ be such a matroid and write $\gamma(N) =: (c_1,\dots,c_r)$. First assume all the $Z_i$ are cyclic in $N$. As $N \neq M$, there must be a minimal $j \geq 1$ and a minimal cyclic flat $G$ of $N$ such that $Z_{j-1} \subsetneq G \subsetneq Z_j$. Hence $c_i = d_i$ for all $i < j$ and $c_j < d_j$, so $\gamma(N) <_\lex \gamma(M)$.
    
    Now let $1 \leq j$ be minimal such that $Z_j$ is not a cyclic flat of $N$. We can assume that $Z_0,\dots,Z_{j-1}$ are all the cyclic flats of $N$ of rank at most $r_{j-1}$, since otherwise we again have $\gamma(N) <_\lex \gamma(M)$. In particular, $c_i = d_i$ for all $i < j$. Since $Z_j$ is not cyclic, we can use a similar argument as in the case $\gamma(M)=(1,\dots,1)$ to see that there must be a cyclic flat $G'$ of $N$ of rank $r_j' < r_j$ such that $Z_{j-1} \subsetneq G' \subsetneq Z_j$. Hence $c_j \leq r_j' - r_{j-1} < r_j - r_{j-1} = d_j$, so $\gamma(N) <_\lex \gamma(M)$.
 \end{proof}
\end{prop}

\subsection{Cyclic reductions}\label{subsection_cyclic_reductions}

This section is dedicated to the notion of cyclic reductions $N$ of a matroid $M$, which are special cases of \emph{rank-preserving weak maps} $N \stackrel{\textnormal{id}}{\to} M$. A weak map between matroids $M,M'$ on ground sets $E,E'$ is a map $\varphi: E \to E'$, such that for all $X \subseteq E$, we have $\rank_M(X) \geq \rank_{M'}(\varphi(X))$. It is rank-preserving if $\rank(M) = \rank(M')$ (see for example \cite[Chapter 7.3]{omatroidtheory} for more on weak maps).

We will mainly be concerned with the question when a flat of $N$ is a flat of the same rank in $M$ and vice versa.

\begin{defn}
  Let $M,N$ be matroids on $E$. We say that $N$ is a \emph{cyclic reduction} of $M$ if $\{\emptyset, \cyc_M(E)\} \subseteq \curly{Z}(N) \subseteq \curly{Z}(M)$ and the rank function on $\curly{Z}(N)$ is the one given by $M$, i.e.\ $\rank_N(Z) = \rank_M(Z)$ for all $Z \in \curly{Z}(N)$.
\end{defn}

\begin{ex}
 The easiest way to create a cyclic reduction of a loopfree matroid $M$ is to pick a chain $\emptyset = Z_0 \subsetneq \dots \subsetneq Z_k = \cyc_M(E)$ and define $N$ to be the matroid with cyclic flats $\curly{Z}(N) = \{Z_i; i=0,\dots,k\}$ with $\rank_N(Z_i) = \rank_M(Z_i)$. In particular, $U_{r,n}$ is a cyclic reduction of any loop- and coloopfree matroid of rank $r$ on $E = [n]$.
\end{ex}

\begin{lemma}\label{lemma_flatcontainer}
 Let $N$ be a cyclic reduction of $M$. For every flat $F$ of $N$, there exists a flat $G$ of $M$ such that $F \subseteq G$, $\rank_N(F) = \rank_M(G)$ and $\frank_M(G) \leq \frank_N(F)$. In particular, $\rank_M(A) \leq \rank_N(A)$ for all $A \subseteq E$, so $N \stackrel{\textnormal{id}}{\to} M$ is a rank-preserving weak map.
 \begin{proof}
  The statement is clearly true if $s := \rank(F) = 0$, i.e.\ $F = \emptyset$. Now assume $s > 0$. If $F$ is cyclic, then we can choose $G = F$. Otherwise, let $m := \frank(F)$. Then by Remark \ref{remark_cyclic_construction} $F$ is of the form $F = F' \cup \{x\}$, where $F' \in \curly{F}_{s-1,m-1}(N)$ and $x \notin F'$. By induction there exists a flat $G' \in \curly{F}(M)_{s-1,j}$ with $j \leq m-1$ and $F' \subseteq G'$. If $x \in G'$, then we can pick any $y \notin G'$. Then there exists a flat $G' \cup \{y\} \subseteq G \in \curly{F}(M)_{s,k}$, with $k \leq j+1 \leq m$. If $x \notin G'$, we can pick $y = x$ and apply the same argument.
 \end{proof}
\end{lemma}

\begin{prop}\label{prop_passdown_flats}
 Let $N$ be a cyclic reduction of $M$. Let $F$ be a flat of $M$. Then $F$ is a flat of $N$ of the same rank if and only if $\cyc_M(F)$ is a cyclic flat of $N$.
 \begin{proof}
  For the \enquote{if} direction, assume $\cyc_M(F)$ is a cyclic flat of $N$. We prove that $F$ is a flat of $N$ by induction on $\frank(F)$. If $\frank(F) = 0$, then $F = \cyc_M(F)$ and we are done. So assume $m := \frank(F) > 0$ and write $s := \rank_M(F)$. Pick any $p \in F \wo \cyc_M(F)$. The set $F' := F \wo \{p\}$ is again a flat of $M$ and $\cyc_M(F') = \cyc_M(F)$. In particular, $F' \in \curly{F}_{s-1,m-1}(M)$, so by induction it is also in $\curly{F}_{s-1,m-1}(N)$. If $F \notin \curly{F}_{s,m}(N)$, then by Remark \ref{remark_cyclic_construction} there must be a flat $G \in \curly{F}_{s,l}(N)$, with $l < m$ and $F \subsetneq G$. By Lemma \ref{lemma_flatcontainer} there exists a flat $H \in \curly{F}_s(M)$ with $F \subsetneq G \subseteq H$, which is clearly impossible. Hence $F \in \curly{F}_{s,m}(N)$.
  
  For the \enquote{only if} direction, assume $Z := \cyc_M(F)$ is not a cyclic flat of $N$, but that $F$ is a flat of both $M$ and $N$, such that $\rank_M(F) = \rank_N(F)$. As $N$ is a cyclic reduction of $M$, we must have $Z' := \cyc_N(F) \subsetneq Z$. We then have
  \begin{align*}
    \rank_M(Z) + \abs{F \wo Z} &= \rank_M(F) = \rank_N(F) \\
    &= \rank_N(Z') + \abs{F \wo Z'} = \rank_M(Z') + \abs{F \wo Z'}\enspace.
  \end{align*}
On the other hand, property (Z2) of Theorem \ref{thm_cyclic_axioms} tell us that 
\begin{align*}
 \abs{F \wo Z'} - \abs{F \wo Z} = \abs{Z \wo Z'} > \rank_M(Z) - \rank_M(Z')\enspace,
\end{align*}
which is a contradiction.
 \end{proof}
\end{prop}

\begin{defn}
 Let $N$ be a cyclic reduction of $M$ and $F \in \curly{F}(N)$ be a flat of $N$. We call $Z \in \curly{Z}(M)$  an \emph{abundant flat for $F$ in $M$} if 
 $$\abs{Z \cap F} \geq \nty_N(F) + \rank_M(Z)\enspace.$$
 Equivalently, $\abs{F \wo Z} \leq \rank_N(F) - \rank_M(Z)$. Note that this implies in particular that $\rank_M(Z) \leq \rank_N(F)$.
 
 Moreover, if $Z \notin \curly{Z}(N)$, we call $Z$ a \emph{witness flat for $F$ in $M$}. We denote the set of all abundant flats for $F$ by $A_M(F)$ and the set of witness flats by $W_M(F)$.
\end{defn}

\begin{remark}\label{remark_cyclic_abundant}
 It follows from the remarks in Definition \ref{def_cyclic_flats} that $\cyc_N(F) \in A_M(F)$. Also, any flat which is strictly contained in $\cyc_N(F)$ can never be an abundant flat for $F$ by axiom (Z2) of Theorem \ref{thm_cyclic_axioms}.
\end{remark}

\begin{prop}\label{prop_liftup_flats}
Let $N$ be a cyclic reduction of $M$. Let $F \in \curly{F}_k(N)$. Then 
$$A_M(F) \wo W_M(F) = \{ \cyc_N(F) \}\enspace.$$ Furthermore, $F \in \curly{F}_k(M)$ if and only if $W_M(F) = \emptyset$.
\begin{proof}
 We know from Remark \ref{remark_cyclic_abundant} that $\cyc_N(F) \in A_M(F)$. So let $Z \in \curly{Z}(N)$ and assume $\abs{Z \cap F} \geq \nty_N(F) + \rank_M(Z) = \abs{F} - (k - \rank_M(Z))$. Using the construction from Remark \ref{remark_cyclic_construction}, we can inductively build flats $G_i$ of $N$ of rank $\rank_M(Z) + i$ such that $\abs{G_i \cap F} \geq \abs{Z \cap F} + i$. In particular, there is a flat $Z \subseteq G$ of rank $k$ such that $F \subseteq G$. Hence $G = F$. By construction of $G$ we have $Z \subseteq F$, so $Z \subseteq \cyc_N(F)$. As we assumed $Z \in A_M(F)$, this implies $\cyc_N(F) = Z$ by Remark \ref{remark_cyclic_abundant}. That proves the first statement. 
 
 For the \enquote{only if} part of the second statement, assume that $F \in \curly{F}_k(M)$ as well. Let $Z \in \curly{Z}(M)$ and assume $Z \in W_M(F)$. We construct a flat $G \in \curly{F}_k(M)$ in the same manner as before and conclude again that $F = G$ and $Z = \cyc_M(F)$. By Proposition \ref{prop_passdown_flats}, $Z \in \curly{Z}(N)$, which is a contradiction.
 
 For the \enquote{if} direction, assume $F \in \curly{F}_k(N) \wo \curly{F}_k(M)$. By Lemma \ref{lemma_flatcontainer}, there exists a flat $F \subsetneq G$ in $M$ of rank $k$. Let $Z := \cyc_M(G)$. Then by Proposition \ref{prop_passdown_flats}, $Z \notin \curly{Z}(N)$, since otherwise $G$ would be a flat of $N$. Let $Z' := \cyc_N(F)$. Then $Z' \subseteq Z$, so  $\abs{Z \cap F} \geq \abs{Z'} = \abs{F} - (k - \rank_M(Z')) \geq \abs{F} - (k - \rank_M(Z)).$ Hence $Z \in W_M(F)$.
 \end{proof} 
\end{prop}

\begin{lemma}\label{lemma_witness_join_or_meet}
Let $N$ be a cyclic reduction of $M$. Let $F \subseteq F'$ be flats of $N$ and assume $Z \in A_M(F)$ and $Z' \in A_M(F')$. Then at least one of the following is true: 
\begin{itemize}
  \item $Z \wedge_{\curly{Z}(M)} Z' \in W_M(F)\enspace.$
  \item $Z \vee_{\curly{Z}(M)} Z' \in A_M(F')\enspace.$
 \end{itemize}
 If $Z' \in W_M(F')$, then either $Z \wedge_{\curly{Z}(M)} Z' \in W_M(F)$ or $Z \vee_{\curly{Z}(M)} Z' \in W_M(F')$.
 \begin{proof}
  As a shorthand, write $X := Z \wedge_{\curly{Z}(M)} Z'$ and $Y := Z \vee_{\curly{Z}(M)} Z'$. Assume that $X \notin W_M(F)$. Then by Proposition \ref{prop_liftup_flats} $\abs{X \cap F} \leq \nty_N(F) + \rank_M(X)$, with equality if and only if $X = \cyc_N(F)$. We have to show that $Y$ contains sufficiently many elements of $F'$. So we compute
 \begin{align*}
  \abs{Y \cap F'} &\geq \abs{ (Z \cup Z') \cap F'} = \abs{Z \cap F'} + \abs{Z' \cap F'} - \abs{Z \cap Z' \cap F'} \\
\textnormal{\tiny as $F \subseteq F'$}\quad&= \abs{Z \cap F} + \abs{Z' \cap F'} - \abs{Z \cap Z' \cap F'} + \abs{Z \cap (F' \wo F)}\\
\textnormal{\tiny as $Z \in A_M(F), Z' \in A_M(F)$}\quad&\geq \nty_N(F) + \nty_N(F') + \rank_M(Z) + \rank_M(Z') \\
		 &\qquad - \abs{Z \cap Z' \cap F'} + \abs{Z \cap (F' \wo F)}\\
 \textnormal{\tiny by (Z3)}\quad&\geq (\nty_N(F) + \rank_M(X)) + (\nty_N(F') + \rank_M(Y))\\
		 &\qquad + \abs{ (Z \cap Z') \wo X} - \abs{Z \cap Z' \cap F'} + \abs{Z \cap (F' \wo F)}\\
\textnormal{\tiny as $X \notin W_M(F)$}\quad&\geq \nty_N(F') + \rank_M(Y)\\
		 &\qquad + \underbrace{\abs{X \cap F}+ \abs{ (Z \cap Z') \wo X} - \abs{Z \cap Z' \cap F'} + \abs{Z \cap (F' \wo F)}}_{=: \delta}\enspace.\\
 \end{align*}
 So we only have to show that $\delta \geq 0$. For this note that it is trivially true that $\abs{Z \cap (F' \wo F)} \geq \abs{X \cap (F' \wo F)}$. Also, we have
 \begin{align*}
\abs{Z \cap Z' \cap F'} - \abs{X \cap F} - \abs{X \cap (F' \wo F)} &= \abs{Z \cap Z' \cap F'} - \abs{X \cap F'} \\
&= \abs{ ((Z \cap Z')\wo X) \cap F'}\enspace.  
 \end{align*}
Hence we have 
 $$\delta \geq \abs{(Z \cap Z') \wo X} - \abs{ ( (Z \cap Z') \wo X) \cap F'} \geq 0\enspace.$$
 Finally, assume that $Z' \in W_M(F')$. As we already proved that $Y \in A_M(F)$, then by Proposition \ref{prop_liftup_flats} we only need to prove that $Y \neq \cyc_N(F')$. But if that was the case, we would have $Z' \subseteq Y = \cyc_N(F')$, so by Remark \ref{remark_cyclic_abundant} $Z'$ cannot be a witness flat for $F'$.
 \end{proof}
\end{lemma}

\subsection{Representations of arbitrary matroids}\label{subsection_representations}

We will begin by stating how an arbitrary matroid can be written as a linear combination of nested matroids. 

\begin{defn}
 Let $M$ be a matroid. For any maximal chain $\curly{C} = (F_0 = \emptyset \subsetneq F_1 \subsetneq \dots \subsetneq F_r = E$ of flats we define its \emph{cyclic set} to be 
 $\cyc(\curly{C}) = \{ \cyc_M(F_i); i = 0,\dots,r\} \subseteq \curly{Z}(M)$. The \emph{cyclic chain lattice} of $M$ is the set 
 $$\TCP(M) := \{ T \subseteq \curly{Z}(M) \textnormal{ a chain with } \emptyset, 1_{\curly{Z}} \in T\} \cup \{ \hat{1}\}\enspace,$$ 
 with partial order induced by set inclusion and $\hat{1}$ as an artificial maximal element, i.e.:
 \begin{itemize}
  \item $T < \hat{1}$ for all chains $T$.
  \item If $T,T' \neq \hat{1}$, then $T \leq T'$ if and only if $T \subseteq T'$.
 \end{itemize}
This is a lattice with $T \wedge T' = T \cap T'$ and 
$$T \vee T' = \begin{cases}
               T \cup T', &\textnormal{if } T \cup T' \textnormal{ is a chain,}\\
               \hat{1}, &\textnormal{otherwise.}
              \end{cases}$$

Note that each element $\hat{1} \neq T \in \TCP(M)$ defines a chain of cyclic flats of $M$ (with ranks given by $\rank_M$) and thus a nested matroid, which we denote by $M\nest{T}$.

Furthermore, we will use the following shorthand: For any $T \in \TCP(M)$, let 
$$\mu_1(T) := \mu_{\TCP(M)}(T,\hat{1})\enspace.$$
\end{defn}

\begin{theorem}\label{theorem_basis_presentation}
 Let $M$ be a loopfree matroid of rank $r$ on $[n]$. Then the following equality holds in $\M_{r,n}$:
 $$M = \sum_{\substack{T \in \TCP(M)\\T \neq \hat{1}}} \left( - \mu_1(T) \right) M\nest{T}\enspace.$$
\begin{proof} 
 Let $\curly{C} = (\emptyset = F_0 \subsetneq F_1 \subsetneq \dots \subsetneq F_r = E)$ be a chain of length $r = \rank(M)$. 
  
  First assume $\curly{C}$ is a chain of flats of $M$ and denote by $T$ its cyclic set. Then by Proposition \ref{prop_passdown_flats}, $\curly{C}$ is a chain of flats in $M\nest{T'}$ for an element $\hat{1} \neq T' \in \TCP(M)$ if and only if $T \leq T'$. So by Lemma \ref{lemma_moebius} the coefficient of $\curly{C}$ on the right hand side in the above equation is
  $$\sum_{\hat{1} \neq T' \geq T} - \mu_1(T') = \mu_1(\hat{1}) = 1\enspace.$$
  
  Now assume $\curly{C}$ is a chain of some $M\nest{T}$, but not of $M$. Let $T_1,\dots,T_k$ be the minimal elements of $\TCP(M)$ such that $\curly{C}$ is a chain in $M\nest{T_i}, i = 1,\dots,k$. In particular, by Proposition \ref{prop_passdown_flats} every $Z \in T_i$ is the cyclic part of some $F_j$. Now if $T$ is such that $\curly{C}$ is a chain in $M\nest{T}$, there is a unique $i = 1,\dots,k$ such that $T \geq T_i$: If $T \geq T_i, T_j$, then by Proposition \ref{prop_passdown_flats}, both $T_i$ and $T_j$ contain $\cyc_{M\nest{T}}(\curly{C})$. Hence so does $T_i \cap T_j$, which is a contradiction to the minimality assumption.
  
  For $i = 1,\dots,k$ we define the following subposet of $\TCP(M)$:
  $$R_i := \{\hat{1} \neq T > T_i;\; \curly{C} \textnormal{ is not a chain in } M\nest{T}\}\enspace.$$
  
  Proposition \ref{prop_passdown_flats} implies that if $T \in R_i$ and $T' > T$, then $T' \in R_i$. It follows from Lemma \ref{lemma_moebius} that 
  \begin{equation}\label{eq_ri_sum}
    \mu(R_i) = \mu_{\hat{R_i}}(\hat{0},\hat{1}) = - \sum_{T \in R_i \cup \{\hat{1}\}} \mu_{\hat{R_i}}(T,\hat{1}) = - \sum_{T \in R_i \cup \{\hat{1}\}} \mu_1(T)\enspace.
  \end{equation}
  
  We write $W_i := \{Z \in \curly{Z}(M); T_i \cup \{Z\} \in R_i\}$. Equivalently, this is the set of cyclic flats $Z \notin T_i$ such that $T_i \cup \{Z\}$ is a chain and $Z \in W_M(F)$ for some $F$ in $\curly{C}$. Note that if $W_i$ is not empty, it is a join-contractible subposet of $\curly{Z}(M)$: Pick any minimal element $Z$ of $W_i$ and let $Z' \in W_i$ be arbitrary. If $Z \subseteq Z'$, then $Z \vee_{\curly{Z}(M)} Z' = Z' \in W_i$. Otherwise, the meet of $Z$ and $Z'$ cannot lie in $W_i$ due to the minimality of $Z$. By Lemma \ref{lemma_witness_join_or_meet} $Z \vee_{\curly{Z}(M)} Z'$ then lies in $W_i$.
  
  By the above considerations and Proposition \ref{prop_liftup_flats}, we have a poset isomorphism
  $$\Ch(W_i) \to R_i,\quad \{Z_1,\dots,Z_k\} \mapsto T_i \cup \{Z_1,\dots,Z_k\}\enspace.$$  
  
  Applying Lemma \ref{lemma_join_contractible} and Proposition \ref{prop_moebius_chains} we see that if $R_i \neq \emptyset$, we have
  \begin{equation}\label{eq_wi_mu}
   \mu(R_i) = \mu(\Ch(W_i)) = \mu(W_i) = 0\enspace.
  \end{equation}

   We want to show that $R_i \neq \emptyset$ or, equivalently, that $W_i \neq \emptyset$. We need to construct a witness flat that forms a chain with $T_i$. As we assumed that $\curly{C}$ is not a chain in $M$, Proposition \ref{prop_liftup_flats} tells us that there must be a noncyclic flat $F$ in $\curly{C}$ and a witness flat $Z \in W_M(F)$. We choose $F$ maximal, such that a witness flat exists for it and let $Z$ also be a maximal element of $W_M(F)$. Denote by $Z_1 := \cyc_{M\nest{T_i}}(F)$. Then $Z_1$ is not equal to $1_{\curly{Z}(M)}$, otherwise there could be no witness flats for $F$. So let $Z_2$ be the smallest element of $T_i$ that strictly contains $Z_1$. As $Z_1 \wedge_{\curly{Z}(M)} Z \subseteq Z_1$, it cannot be a witness flat for $F$ by Remark \ref{remark_cyclic_abundant}. Hence, by Lemma \ref{lemma_witness_join_or_meet} we know that $Z_1 \vee_{\curly{Z}(M)} Z \in W_M(F)$. Due to the maximality of $Z$, we must have $Z_1 \vee_{\curly{Z}(M)} Z = Z$, so $Z_1 \subseteq Z$. Due to the minimality of $T_i$, there must be a flat $F \subsetneq F'$ such that $Z_2 = \cyc_{M\nest{T_i}}(F')$. In particular, $Z_2 \in A_M(F')$. 
   Let $X := Z_2 \wedge_{\curly{Z}(M)} Z$ and $Y := Z_2 \vee_{\curly{Z}(M)} Z$. By Lemma \ref{lemma_witness_join_or_meet} there are two possibilities: If $X \in W_M(F)$, we have $Z_1 \subseteq X \subseteq Z_2$, so $X \in W_i$. If $Y \in A_M(F')$, then by our choice of $F$ there are no witness flats for $F'$, so $Y = Z_2$ by Proposition \ref{prop_liftup_flats}. But then $Z_1 \subseteq Z \subseteq Y = Z_2$, so $Z \in W_i$.
  
   We can now finally use Equations \ref{eq_ri_sum} and \ref{eq_wi_mu}, as well as Lemma \ref{lemma_moebius} to see that the coefficient of $\curly{C}$ on the right hand side above is
  $$ - \sum_{i=1}^k\left( \sum_{T \geq T_i} \mu_1(T) - \sum_{T' \in R_i \cup \{\hat{1}\}} \mu_1(T') \right) = - \sum_{i=1}^k (0 - 0) =0\enspace.$$
\end{proof}
 \end{theorem}

\begin{corollary}\label{corollary_basis}
 The set of loopfree nested matroids of rank $r$ on $[n]$ is a basis for $\M_{r,n}$, i.e.\ every matroid can be written as a unique linear combination of nested matroids.
\end{corollary}

\begin{remark}
 It should be noted that Derksen and Fink show in \cite{dfvaluative} that the polytopes of nested matroids (which they call \emph{Schubert matroids}) also form a basis for their module of polytopes. In their case the representation of an arbitrary matroid is given as a sum over \emph{all} possible chains of sets (for details see their Theorem 4.2). It would be interesting to study the precise relation of these two presentations and what it implies for the combinatorics of the matroids involved (see also Remark  \ref{remark_derksen_fink}).
\end{remark}

\begin{ex}\label{ex_representation}\indent\par
\begin{enumerate}
 \item We recall the matroid $M$ on $E = \{1,\dots,4\}$ from example \ref{ex_cyclic_flats}, whose flats were given by $\curly{F}(M) = \curly{Z}(M) = \{\emptyset, \{1,4\}, \{2,3\},E\}$. The theorem tells us that
 $$M = M\nest{\emptyset, \{1,4\}, E} + M\nest{\emptyset, \{2,3\},E} - M\nest{\emptyset,E}\enspace,$$
 which is the same relation we already encountered in Example \ref{ex_relation}.
 \item We consider the matroid $M$ of rank 4 on $E := \{1,\dots,8\}$ given by the lattice of cyclic flats depicted in Figure \ref{figure_torsion_poset}. The cyclic chain poset is drawn below. According to Theorem \ref{theorem_basis_presentation}, we have
 \begin{align*}
M = &\;M\nest{S_1,R,U_1} + M\nest{S_2,R,U_1} + M\nest{S_1,R,U_2} + M\nest{S_2,R,U_2}\\ &\; - M\nest{R,U_1} - M\nest{S_1,R} - M\nest{S_2,R} - M\nest{R,U_2} \\&\;+ M\nest{R}  
 \end{align*}
(Note that we omit $\emptyset,E$ in the description of each chain set).
 
 \begin{figure}[ht]
  \centering
  \begin{tikzpicture}
  \newcommand*{\stepup}{-1.7}
  \newcommand*{\stepleft}{-2}
   \matrix[row sep=20pt]{
      \draw (0,0) -- (-1,1) -- (0,2) -- (-1,3) -- (0,4);
      \draw (0,0) -- (1,1) -- (0,2) -- (1,3) -- (0,4);
      \fill[black] (0,0) circle (2pt) node[below]{\tiny $\emptyset$};
      \fill[black] (-1,1) circle (2pt) node[left]{\tiny $S_1 = \{1,2\}, r(S_1) = 1$};
      \fill[black] (0,2) circle (2pt) node[left]{\tiny $R = \{1,2,3,4\}, r(R) = 2$};
      \fill[black] (-1,3) circle (2pt) node[left]{\tiny $U_1 = \{1,2,3,4,5,6\}, r(U_1) = 3$};
      \fill[black] (0,4) circle (2pt) node[above]{\tiny $E, r(E) = 4$};
      \fill[black] (1,1) circle (2pt) node[right]{\tiny $S_2 = \{3,4\}, r(S_2) = 1$};
      \fill[black] (1,3) circle (2pt) node[right]{\tiny $U_2 = \{1,2,3,4,7,8\}, r(U_2) = 3$}; \\
      \node[rounded corners =2pt,draw] (zero) at (0,0) {\tiny $\substack{\hat{1}\\1}$};
      \node[rounded corners =2pt,draw] (m1) at (-1.5*\stepleft,\stepup) {\tiny $\substack{\{S_1,R,U_1\}\\-1}$};
      \node[rounded corners =2pt,draw] (m2) at (-.5*\stepleft,\stepup) {\tiny $\substack{\{S_2,R,U_1\}\\-1}$};
      \node[rounded corners =2pt,draw] (m3) at (.5*\stepleft,\stepup) {\tiny $\substack{\{S_1,R,U_2\}\\-1}$};
      \node[rounded corners =2pt,draw] (m4) at (1.5*\stepleft,\stepup) {\tiny $\substack{\{S_2,R,U_2\}\\-1}$};
      \node[rounded corners =2pt,draw] (c1) at (-2.4*\stepleft,2*\stepup) {\tiny $\substack{\{S_1,U_1\}\\0}$};
      \node[rounded corners =2pt,draw] (c2) at (-1.7*\stepleft,2*\stepup) {\tiny $\substack{\{S_1,U_2\}\\0}$};
      \node[rounded corners =2pt,draw] (c3) at (-1*\stepleft,2*\stepup) {\tiny $\substack{\{R,U_1\}\\1}$};
      \node[rounded corners =2pt,draw] (c4) at (-.3*\stepleft,2*\stepup) {\tiny $\substack{\{S_1,R\}\\1}$};
      \node[rounded corners =2pt,draw] (c5) at (.3*\stepleft,2*\stepup) {\tiny $\substack{\{S_2,R\}\\1}$};
      \node[rounded corners =2pt,draw] (c6) at (1*\stepleft,2*\stepup) {\tiny $\substack{\{R,U_2\}\\1}$};
      \node[rounded corners =2pt,draw] (c7) at (1.7*\stepleft,2*\stepup) {\tiny $\substack{\{S_2,U_1\}\\0}$};
      \node[rounded corners =2pt,draw] (c8) at (2.4*\stepleft,2*\stepup) {\tiny $\substack{\{S_2,U_2\}\\0}$};
      \node[rounded corners =2pt,draw] (t) at (0,3*\stepup) {\tiny $\substack{\{R\}\\-1}$};
      \node[rounded corners =2pt,draw] (s1) at (-2*\stepleft,3*\stepup) {\tiny $\substack{\{S_1\}\\0}$};
      \node[rounded corners =2pt,draw] (u1) at (-1*\stepleft,3*\stepup) {\tiny $\substack{\{U_1\}\\0}$};
      \node[rounded corners =2pt,draw] (u2) at (1*\stepleft,3*\stepup) {\tiny $\substack{\{U_2\}\\0}$};
      \node[rounded corners =2pt,draw] (s2) at (2*\stepleft,3*\stepup) {\tiny $\substack{\{S_2\}\\0}$};
      \node[rounded corners =2pt,draw] (u) at (0,4*\stepup) {\tiny $\substack{-\\0}$};
      \foreach \x in {(m1), (m2),(m3),(m4)}
      {\draw (zero) -- \x;}
      \foreach \x in {(c3),(c4),(c5),(c6)}
      {\draw (t) -- \x;}
      \foreach \x in {(c1),(c2),(c4)}
      {\draw (s1) -- \x;}
      \foreach \x in {(c1),(c3),(c7)}
      {\draw (u1) -- \x;}
      \foreach \x in {(c2),(c6),(c8)}
      {\draw (u2) -- \x;}
      \foreach \x in {(c5),(c7),(c8)}
      {\draw (s2) -- \x;}
      \foreach \x in {(s1),(s2),(u1),(u2),(t)}
      {\draw (u) -- \x;}
      \draw (c1) -- (m1);
      \draw (c2) -- (m3);
      \draw (c3) -- (m1); 
      \draw (c3) -- (m2);
      \draw (c4) -- (m1);
      \draw (c4) -- (m3);
      \draw (c5) -- (m2);
      \draw (c5) -- (m4);
      \draw (c6) -- (m3);
      \draw (c6) -- (m4);
      \draw (c7) -- (m2);
      \draw (c8) -- (m4);
      \\
   };
  \end{tikzpicture}
  \caption{The lattice of cyclic flats of a matroid and the corresponding cyclic chain lattice with the value of the Möbius function $\mu_1(\cdot)$ indicated beneath each element. Note that each cyclic set naturally contains $\emptyset,E$ but we leave them out here for legibility.}\label{figure_torsion_poset}
 \end{figure}
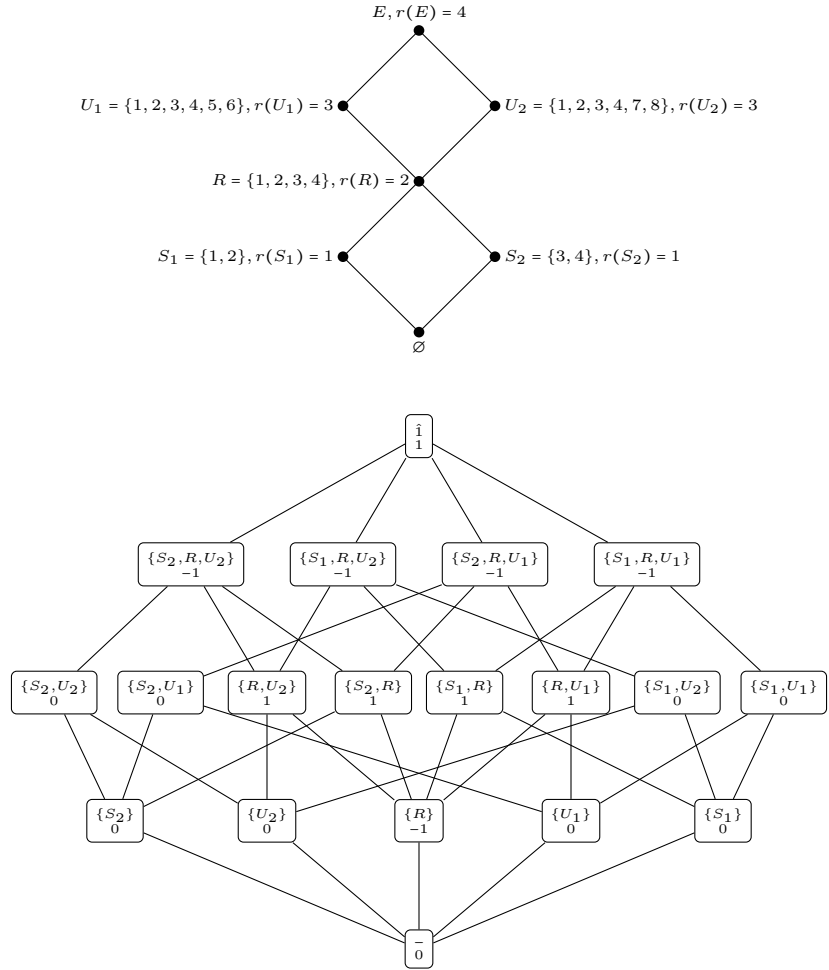
\end{enumerate}
\end{ex}

\subsection{The \texorpdfstring{$\curly{G}$}{G}-invariant}\label{subsection_tutte}
%
%
%

As mentioned in the introduction, the $\mathcal{G}$-invariant is an important matroid invariant. We will use the characterization by Bonin and Kung \cite{bkcatenary} in terms of \emph{catenary data} as a definition. 

\begin{defn} Assume $n \geq 1$ and $1 \leq r \leq n$. Let $\curly{C} = (F_0 \subsetneq F_1 \subsetneq \dots \subsetneq F_r = E)$ be a chain of sets. Its \emph{composition} is the tuple $(a_0,\dots,a_r)$ with $a_0 := \abs{F_0}$ and $a_i := \abs{F_i \wo F_{i-1}}$ for $i > 0$. For a matroid $M$ of rank $r$ and a fixed composition $a := (a_0,\dots,a_r)$ we define $\nu(M;a)$ to be the number of chains of flats of $M$ of length $r$ with composition $a$. 

Now let $\curly{G}(n,r)$ be the free abelian group on all possible compositions $a$. We will denote the generator corresponding to $a$ by $\gamma(a)$ in accordance with the notation in \cite{bkcatenary}. The \emph{$\curly{G}$-invariant} of $M$ is 
$$\curly{G}(M) := \sum_a \nu(M;a) \gamma(a) \in \curly{G}(n,r)\enspace,$$
where the sum runs over all possible compositions $a$.
\end{defn}

\begin{theorem}\label{theorem_g_invariant}
 For each $n$ and $r$ the $\curly{G}$-invariant induces a $\Z$-module homomorphism via
 $$\M_{r,n} \to \curly{G}(r,n),\; M \mapsto \curly{G}(M)\enspace.$$
 \begin{proof}
  We only need to show that for each composition $a$, there is an induced $\Z$-module homomorphism $\M_{r,n} \to \Z,\, M \mapsto \nu(M;a)$.  But this is obvious from the definition of $\M_{r,n}$ as a submodule of $V_{r,n}$: Project to the coordinates of chains with composition $a$, then take the sum of the coordinates.
 \end{proof}
\end{theorem}

The significance of the $\mathcal{G}$-invariant can also be recognized from the fact that many other matroid invariants can be derived from it. In particular, Derksen \cite{dsymmetricpolymatroids} showed that the \emph{Tutte polynomial} 
$$t_M(x,y) := \sum_{S \subseteq E} (x-1)^{\corank_M(S)} (y-1)^{\nty_M(S)} \in \Z[x,y]$$
can be computed from the $\mathcal{G}$-invariant by applying a certain linear map to $\mathcal{G}(n,r)$. This immediately implies the following:

\begin{corollary}
 The Tutte polynomial induces a $\Z$-module homomorphism via
 $$\M_{r,n} \to \Z[x,y],\; M \mapsto t_M(x,y)\enspace.$$
\end{corollary}

\begin{remark}
This implies of course that the same statement holds for all \emph{generalized Tutte-Grothendieck invariants} and \emph{Tutte-Grothendieck group invariants}, as defined in \cite{botutteapplications}, such as for example the characteristic polynomial, the beta invariant, the Whitney numbers of the first kind or the number of bases.

However, there are also invariants which can not be derived from the Tutte polynomial, but which are encoded in the $\curly{G}$-invariant. We consider a particular example from \cite{bkcatenary}, again adopting their notation:
\end{remark}

\begin{defn}
  For a matroid $M$, we denote by $F_{h,k}(M;s_h,\dots,s_k)$ the number of chains of flats $F_h \subsetneq \dots \subsetneq F_k$ of $M$ such that $\rank(F_i) = i$ and $\abs{F_i} = s_i$ for all $i = h,\dots,k$.
\end{defn}

\begin{lemma}
 For fixed $h,k$ and $s := (s_h,\dots,s_k)$, there is a $\Z$-module homomorphism $\M_{r,n} \to \Z$ induced by $M \mapsto F_{h,k}(M;s_h,\dots,s_k)$.
 \begin{proof}
  In \cite[Proposition 5.2]{bkcatenary}, Bonin and Kung show how the numbers $F_{h,k}(M;s)$ are derived from the $\mathcal{G}$-invariant. We follow their argument to see that we indeed obtain a linear map.
  
  For $k = r$, they show there exists a linear \emph{specialization map} $\textnormal{spec}: \mathcal{G}(n,r) \to \Z$ such that the following identity holds: $$F_{h,r}(M;s) = \textnormal{spec}\left((1/s_h!)\sum_a \nu(M; a') \gamma(a)\right)\enspace,$$
  where $a$ runs over a certain set of compositions depending on $h$ and $s$ and $a'$ is a composition that depends only on $a$ and $s$. We saw in Theorem \ref{theorem_g_invariant} that $\nu(M;a')$ induces a linear map on $\M_{n,r}$. Hence $F_{h,r}(M;s)$ induces a linear map.
  
  For $k = r-1$ the statement follows from the observation that $$F_{h,r-1}(M;s) = F_{h,r}(M;s,n)\enspace.$$
  
  For $k < r-1$, note that $F_{h,k}(M;s) = F_{h,k}(T^{r-k-1}M;s)$, where $T^iM$ denotes again the $i$-fold truncation of $M$. Recall from example \ref{ex_intersection} that $T^iM = M \cdot U_{n-i,n}$, so $M \mapsto T^iM$ induces a linear map. Thus the claim follows by induction.
 \end{proof}
\end{lemma}

\begin{corollary}
 Let $f_k(M;s,c)$ denote the number of flats $F$ of $M$ of size $s$ and rank $k$ such that $M_{\mid F}$ has $c$ coloops. This induces a $\Z$-module homomorphism $\M_{r,n} \to \Z$. In particular, the number of flats of rank $k$ and the number of cyclic flats of rank $k$ induce $\Z$-module homomorphisms as well.
 \begin{proof}
  The statement is trivial for $k = r$. It is shown in \cite[Proposition 5.5]{bkcatenary} that for $k=r-1$ we have
  $$\sum_{j=c}^n f_{r-1}(M;s,j) \frac{j!}{(j-c)!} = F_{r-1-c,r-1}(M;s-c,s-c+1,\dots,s)\enspace.$$
  An easy inductive argument then shows that for $k=r-1$ we obtain a linear map. The general case then follows by using truncations as in the proof above.
 \end{proof}
\end{corollary}

\begin{ex}
 We again look at the matroid $M$ on $E = \{1,\dots,4\}$ from Examples \ref{ex_cyclic_flats} and \ref{ex_representation},(1). From the latter we recall the linear relation
 $$M = M\nest{\emptyset,\{1,4\},E} + M\nest{\emptyset,\{2,3\},E} - M\nest{\emptyset,E}\enspace.$$
 Computing the $\curly{G}$-invariants on the right hand side, we see that 
 \begin{align*}
    \curly{G}(M) &= ( 2\gamma(0,1,3) + \gamma(0,2,2)) + (2\gamma(0,1,3) + \gamma(0,2,2)) - 4 \gamma(0,1,3)\\
        &= 2\gamma(0,2,2)\enspace.
 \end{align*}
 Indeed, $M$ has two maximal chains of flats, $(\emptyset, \{1,4\},E)$ and $(\emptyset,\{2,3\},E)$.
\end{ex}

\section{Counting nested matroids}\label{section_counting_nested}

We have shown that nested matroids are a basis for $\M_{r,n}$, so naturally we want to determine their exact number. It was already pointed out in \cite{oprmatroidsdomains} that the number of \emph{isomorphism classes} of nested matroids of rank $r$ on $n$ elements is $\binom{n}{r}$. We will show that without taking isomorphisms into account, we still get a familiar number.

\begin{defn}
 Let $n \geq 1$ and $0 \leq r < n$. The \emph{Eulerian number} $A_{r,n}$ is the number of permutations on $\{1,\dots,n\}$ with $r$ \emph{ascents}. An ascent of a permutation $\sigma$ is a number $i \in \{1,\dots,n-1\}$, such that $\sigma(i) < \sigma(i+1)$.
\end{defn}

\begin{remark}
 Another interpretation of these numbers is the following: 
 
 For $0 < k < n$, define the \emph{hypersimplex}
 $$\Delta_{k,n} := \left\{(x_1,\dots,x_n), x_i \in [0,1], \sum x_i  = k\right\}\enspace.$$
 Then $A_{r,n}$ is the lattice volume of $\Delta_{r+1,n+1}$ (see \cite{lpalcoved} for a discussion of a concrete unimodular triangulation of the hypersimplex).
 
 We note a few properties of these numbers here, which can be found in most standard textbooks on combinatorics such as \cite{gkpconcretemath}:
 \begin{itemize}
  \item By definition, we have $\sum_{r=0}^{n-1} A_{r,n} = n!$.
  \item The Eulerian numbers are symmetric: $A_{r,n} = A_{n-r-1,n}$.
  \item The generating function of the Eulerian numbers is
  $$\sum_{n=0}^\infty \sum_{r=0}^\infty A_{r,n}  \frac{x^rt^n}{n!} = \frac{x-1}{x - e^{(x-1)t}}\enspace.$$
 \end{itemize}
 \end{remark}
 
 \begin{defn}
  We denote by $N_{r,n}$ the number of loopfree nested matroids of rank $r$ on $n$ labeled elements.
 \end{defn}
 
 \begin{lemma}\label{lemma_nested_recursive}
 $N_{r,n}$ is determined by the following recursive relation:
 $$N_{r,n} = 1 + \sum_{k=1}^{r-1} \left(\sum_{s=k+1}^{k+n-r} \binom{n}{s} N_{r-k,n-s}\right) \textnormal{ if }r > 1\enspace,$$
 and $N_{1,m} = 1$ for any $m \geq 1$.
 \begin{proof}
  There is only one loopfree matroid of rank 1 on $m$ elements, the uniform matroid, which is also nested. Hence $N_{1,m} = 1$.
  
  We build a nested matroid of rank $r > 1$ by recursively constructing its chain of cyclic flats. The first one is always the empty set of rank zero, as the matroid is loopfree. We consider the choices we have for the first nonempty cyclic flat $F$: We can choose its rank $k$ and its size $s$. If the rank is $r$, then $F = E$ and there is only one nested matroid of this form. Otherwise pick any rank $1 \leq k < r$. By (Z2) of Theorem \ref{thm_cyclic_axioms}, the size $s$ of $F$ has to be at least $k+1$ and can be at most $k+n-r$. We clearly have $\binom{n}{s}$ possibilities to choose such an $F$. For each nested matroid of rank $r$ on $[n]$ whose chain of cyclic flats is of the form$(\emptyset,0) \subsetneq (F_1 = F, r_1 = k) \subsetneq \dots \subsetneq (F_l,r_l)$, the deletion of $F$ gives a nested matroid of rank $r-k$ on $n-s$ elements. Conversely, we can take any such matroid and lift it to a nested matroid of rank $r$ on $n$ elements, whose chain of cyclic flats starts with $(\emptyset,0) \subsetneq (F,k)$. This proves the claim.
 \end{proof}
\end{lemma}

\begin{theorem}\label{theorem_counting_nested}
 For any $n \geq 1$ and $1 \leq r \leq n$, we have
 $$N_{r,n} = A_{r-1,n}\enspace.$$
 \begin{proof}
  We prove this by showing equality of generating functions. As stated above, we have
  $$\sum_{n=0}^\infty \sum_{r=0}^\infty A_{r,n}  \frac{x^rt^n}{n!} = \frac{x-1}{x - e^{(x-1)t}}\enspace.$$
  Here we set $A_{0,0} = 1$ and $A_{r,n} = 0$ if $0 < n \leq r$ or $\min\{n,r\} < 0$.
  
  In accordance with this, we define 
  $$N_{1,0} = 1\textnormal{ and }N_{a,b} = 0\textnormal{ for all }(a,b) \notin \{(r,n); 1 \leq r \leq n, 1 \leq n\} \cup \{(1,0)\}\enspace.$$ We can then rewrite the formula of Lemma \ref{lemma_nested_recursive} as
  $$N_{r,n} = \sum_{k=1}^\infty \sum_{s=k+1}^\infty \binom{n}{s} N_{r-k,n-s}, \textnormal{ if } r > 1\enspace.$$
We thus compute:
\begin{align*}
 f(x,t) &= \sum_{n=0}^\infty \sum_{r=0}^\infty N_{r+1,n} \frac{x^r t^n}{n!}\\
 &= \sum_{n=0}^\infty N_{1,n} \frac{t^n}{n!} + \sum_{n=0}^\infty \sum_{r=1}^\infty \left( \sum_{k=1}^\infty \sum_{s=k+1}^\infty \binom{n}{s} N_{r-k+1,n-s}\right) \frac{x^rt^n}{n!}\\
 &= e^t + \sum_{k=1}^\infty \sum_{s=k+1}^\infty  \frac{x^kt^s}{s!}\underbrace{\left( \sum_{n=0}^\infty \sum_{r=1}^\infty N_{r-k+1,n-s} \frac{x^{r-k} t^{n-s}}{(n-s)!} \right)}_{= f(x,t)}\\
 &= e^t + f(x,t)\left( \sum_{k=1}^\infty \left( e^t - \sum_{i=0}^k \frac{t^i}{i!} \right)x^k\right) \\
 &= e^t + f(x,t)\left(\frac{xe^t - e^{tx}}{1-x} + 1\right)\enspace.\\
\end{align*}
Solving for $f$, we get
$$f(x,t) = \frac{-e^t (1-x)}{xe^t - e^{tx}} = \frac{x-1}{x - e^{(x-1)t}}\enspace.$$
 \end{proof}
\end{theorem}

\section{The intersection product}\label{section_intersection_product}

In this section, we will study the properties of the intersection product on $\M_n$. We will first prove that each nested matroid is a certain product of corank one matroids. We will then study when certain products of matroids vanish in $\M_n$. Finally we will show that $\M_n$ fulfills a Poincar\'e-type duality.

\subsection{Chain products}\label{subsection_chain_matroids}

\begin{defn}
 Fix a ground set $E$ of size $n$. One sees easily that a loopfree corank one matroid is uniquely determined by fixing its set of coloops $G$, where $0 \leq \abs{G} \leq n-2$. It must then be of the form
 $$H_G := U_{\abs{G},G} \oplus U_{\abs{G^c}-1,G^c}\enspace.$$
 Given a chain of sets $\curly{G} = (G_1 \subsetneq \dots \subsetneq G_k)$ with $\abs{G_k} \leq n-2$, we define its \emph{chain product} to be
 $$M_{\curly{G}} := H_{G_1} \wedge \dots \wedge H_{G_k}\enspace.$$
 We also fix the following notation for set systems: For a set $S$ and a positive integer $l$, we write $S^{\oplus l}$ for the $l$-fold concatenation of $S$, i.e. 
 $$S^{\oplus l} = \underbrace{(S,\cdots,S)}_{l \textnormal{ times}}\enspace.$$
\end{defn}

\begin{remark}
 We make a few observations about chain products and corank one matroids that are easily verified:
 \begin{itemize}
  \item $M_{\curly{G}}$ is a loopfree matroid of rank $n-k$. In fact, its bases are given by 
  $$\curly{B}(M_{\curly{G}}) = \left\{E\wo \{j_1,\dots,j_k\}; j_i \notin G_i \textnormal{ for all }i \textnormal{ and }\abs{\{j_1,\dots,j_k\}} = k\right\}.$$
  Hence we have $M_{\curly{G}} = H_{G_1} \cdot \,\dots\, \cdot H_{G_k}$.
  \item For any loopfree matroid $M$ with $\rank(M) > 1$ and any corank one matroid $H_G$, their matroid intersection fulfills $\rank(M \wedge H_G) = \rank(M) -1$ (it may have loops, though).
 \end{itemize}
\end{remark}

\begin{lemma}\label{lemma_chain_bases}
 The bases of $M_{\curly{G}}$ are
 $$\left\{B \in \binom{[n]}{n-k}; \abs{G_i \wo B} \leq i - 1 \textnormal{ for all } i = 1,\dots,k\right\}\enspace.$$
 \begin{proof}
  By our remark above the bases of $M_{\curly{G}}$ are of the form $E \wo \{j_1,\dots,j_k\}$, with $j_i \in G_i^c$. In particular, $G_i \wo B$ can contain at most $j_1,\dots,j_{i-1}$. Thus any basis of $M_{\curly{G}}$ is of the given form. 
  
  Conversely let $B \in \binom{[n]}{n-k}$ with $\abs{B \cap G_i} \geq \abs{G_i} - (i - 1)$ for all $i$. In particular $B = E \wo \{j_1,\dots,j_k\}$ for some $j_i \in [n]$. We define an ordering on $[n]$ in the following manner. For $j \in [n]$, let $m_{\curly{G}}(j) := \min\{i; j \in G_i\}$. Then we say that $j <_{\curly{G}} j'$ if and only if $m_{\curly{G}}(j) < m_{\curly{G}}(j')$ or equality holds and $j < j'$. We can assume without loss of generality that $j_1 <_{\curly{G}} \dots <_{\curly{G}} j_k$. By assumption, $j_1$ must lie in $G_1^c$, so $m_{\curly{G}}(j_1) \geq 2$. Using the fact that $\abs{B^c \cap G_2} \leq 1$, one sees that $j_2$ lies in $G_2^c$. One can continue inductively to see that in fact $j_i \in G_i^c$ for all $i$, which concludes the proof.
 \end{proof}
\end{lemma}

\begin{prop}\label{prop_matroid_transversal}
 Let $\curly{G} = (G_1 \subsetneq \cdots \subsetneq G_k)$ with $\abs{G_k} \leq n-2$ and $k \geq 1$. Then $M_{\curly{G}}$ is a transversal matroid. More precisely:
 $$M_{\curly{G}} = M[\curly{A}_{\curly{G}}] := M\left[ G_1^{\oplus \abs{G_1}}, G_2^{\oplus \abs{G_2 \wo G_1}-1},\cdots, G_k^{\oplus \abs{G_k \wo G_{k-1}}-1}, E^{\oplus \abs{G_k^c}-1}\right]\enspace.$$
\begin{proof}
For this one only needs to verify that the bases of $M[\curly{A}_{\curly{G}}]$ are of the form given in Lemma \ref{lemma_chain_bases}, which is obvious.
\end{proof}
\end{prop}

\begin{remark}\label{remark_cardinality_of_coloops}
 Pick any chain $\curly{G} = (G_1,\dots,G_k)$ with $\abs{G_i} = i-1.$ Then the proposition above tells us that
 $$M_{\curly{G}} = M\left[E^{\oplus n-k} \right] = U_{n-k,n}\enspace.$$
 So we see that one matroid can have multiple representations as a chain product. However, the transversal presentation given above identifies it uniquely, as we will shortly prove. This tells us that for any two representations
 $$N = H_{G_1} \cdot \dots \cdot H_{G_k} = H_{G_1'} \cdot \dots \cdot H_{G_k'}\enspace,$$
 we must have $\abs{G_i} = \abs{G_i'}$ for all $i = 1,\dots,k$. Also if $\abs{G_i \wo G_{i-1}} > 1$, we actually must have $G_i = G_i'$.
\end{remark}

\begin{lemma}\label{lemma_maximal_presentation}
 The presentation $\curly{A}_{\curly{G}}$ is maximal.
 \begin{proof}
  By Remark \ref{remark_transversal_presentations}, maximality is fulfilled if we can show that for all $i = 1,\dots,k$, the matroid $M_i := M[\curly{A}_{\curly{G}}]_{\mid G_i^c}$ is coloop-free.
  
  We can rewrite Lemma \ref{lemma_chain_bases} to see that the bases of $M_{\curly{G}}$ are
  $$\left\{B \in \binom{[n]}{n-k}; \abs{B \cap G_i^c} \leq \abs{G_i^c} -k + (i - 1) \textnormal{ for all } i = 1,\dots,k\right\}\enspace.$$
  In addition, one can easily construct a basis where the above inequality is an equality (i.e., equivalently, $\abs{G_i \wo B} = i-1$ for all $i$). Hence the bases of $M_i$ are all subsets of $G_i^c$ of size $\abs{G_i^c} -k + (i-1)$, so $M_i$ is a uniform matroid of corank $k-i+1 \geq 1$, which is coloop-free.
 \end{proof}
\end{lemma}

\begin{prop}\label{prop_chain_matroids_are_nested}
 The set of chain products is equal to the set of loopfree nested matroids.
 \begin{proof}
  From Proposition \ref{prop_matroid_transversal} and Theorem \ref{theorem_nested_matroids} it is obvious that any chain product is a nested matroid. Now let $M := M[A_1,\dots,A_k]$ be a nested matroid and assume without loss of generality that $k = \rank(M)$. Since $M$ is loopfree, we must have $A_k = E$. Assume there is a $j$ such that $A_{j+1} \wo A_j = \{x\}$ for some $x \in E$. Then $M_{\mid A_j^c} = M[ \{x\}, A_{j+2}\wo A_j,\dots,A_k \wo A_j]$ has $x$ as a coloop so we can augment the presentation and replace $A_j$ by $A_j \cup \{x\}$ without changing the matroid. Using similar arguments, one can finally assume that $M = M[A_1^{\oplus t_1},\dots,A_k^{\oplus t_k}]$, such that
  \begin{itemize}
   \item $\emptyset \subsetneq A_1 \subsetneq \cdots \subsetneq A_k = E$.
   \item $\abs{A_i \wo A_{i-1}} \geq 2$ for all $i > 1$.
   \item $0 < t_1 \leq \abs{A_1}$.
   \item $0 < t_i < \abs{A_i \wo A_{i-1}}$ for all $i > 1$.
   \item $\sum_{i=1}^k t_i = \rank(M)$.
  \end{itemize}
  One can then easily construct a chain that produces a chain product with this presentation.
 \end{proof}
\end{prop}

\begin{corollary}\label{corollary_corank_one}
 The ring $\M_n$ is generated in corank one, more precisely: every matroid can be written as a linear combination of products of corank one matroids.
 \begin{proof}
  By Theorem \ref{theorem_basis_presentation}, every matroid is a linear combination of nested matroids. By Proposition \ref{prop_chain_matroids_are_nested}, every nested matroid is a product of corank one matroids.
 \end{proof}
\end{corollary}

\subsection{Vanishing conditions}\label{subsection_vanishing}

In this section we study when certain products in $\M_n$ vanish. More precisely, we will give necessary and sufficient criteria for a product of a matroid and a nested matroid (written as a chain product) to be zero.

\begin{lemma}\label{lemma_vanishing}
 Let $M$ be a matroid on $E$ of rank at least 2. Let $H_G$ be a corank one matroid. Then the following are equivalent:
 \begin{enumerate}
  \item $M \cdot H_G = 0$.
  \item $M$ has a rank one flat $F$ such that $F \cup G = E$.
 \end{enumerate}
\begin{proof}
 Assume $M$ has a rank one flat $F$ such that $F \cup G = E$. Pick any element $f \in F$. As the elements from $F$ are parallel, any basis $B$ of $M$ containing $f$ must fulfill $B \cap F = \{f\}$. In other words, $B \wo \{f\} \subseteq G$. Since any basis $B'$ of $H_G$ containing $f$ also contains all of $G$, we have $B \cap B' = B$, which is too large to be a basis of $M \wedge H_G$. Hence $f$ is a loop of $M \wedge H_G$ and $M \cdot H_G = 0$.
 
 Conversely, assume $M$ has no such flat and pick $e \in E$. Then there exists an element $e' \notin G$, such that $\{e,e'\}$ is independent in $M$. Thus there is a basis $B$ of $M$ containing $\{e,e'\}$. $B' := E \wo \{e'\}$ is a basis of $H_G$, so $B \cap B' = B \wo  {e'}$ is a basis of $M \wedge H_G$ containing $e$. Hence, $M \wedge H_G$ is loopfree.
\end{proof}
\end{lemma}

\begin{lemma}\label{lemma_flats_of_product}
 Let $M$ be a matroid of rank at least 2 and $H_G$ a corank one matroid such that $M \cdot H_G \neq 0$. Then the flats of the intersection product are given by:
 $$\curly{F}(M \cdot H_G) = \{F \in \curly{F}(M) \cap \curly{F}(H_G); M/F \cdot H_G/F \neq 0\}.$$
 An equivalent formulation is the following: Let $F$ be a flat of both $M$ and $H_G$. Then $F$ is also a flat of $M \cdot H_G$ if and only if one of the following two conditions is met:
 \begin{enumerate}
  \item $F \cup G = E$.
  \item There is no flat $F'$ of $M$ that covers $F$ and such that $F' \cup G = E$.
 \end{enumerate}
 \begin{proof}
  As $M' := M \cdot H_G = M \wedge H_G$ is loopfree, a set $F$ is a flat of $M'$ if and only if $M'/F$ is loopfree. But $M'/F = (M \wedge H_G)/F = (M/F) \wedge (H_G/F)$. As $F$ is a flat of both $M$ and $H_G$, both factors are loopfree. This proves the first statement.
  
  The second statement follows from the first, Lemma \ref{lemma_vanishing} and the fact that
  $$H_G / F = \begin{cases}
               U_{\abs{G \wo F}, G \wo F},&\textnormal{if } F \cup G = E,\\
               H_{G \wo F}, &\textnormal{otherwise.}
              \end{cases}$$
  If $F \cup G = E$, then $M/F \cdot H_G/F = M/F \neq 0$, so $F$ is a flat of $M \cdot H_G$. Otherwise, $F$ is a flat if and only if $M/F \cdot H_{G \wo F} \neq 0$. As the rank one flats of $M/F$ are all flats $F' \wo F$, where $F'$ covers $F$ in $\curly{F}(M)$, this is equivalent to the second condition above. 
 \end{proof}
\end{lemma}

\begin{lemma}\label{lemma_rank_down}
 Let $M$ be a matroid on $[n]$, $H_G$ a corank one matroid and assume $M \cdot H_G \neq 0$. Let $\emptyset \neq F$ be a flat of both $M$ and $M \cdot H_G$. Then 
 $$\rank_{M \cdot H_G}(F) = \begin{cases}
                             \rank_M(F) -1, &\textnormal{if } F \cup G = E\\
                             \rank_M(F), &\textnormal{otherwise.}
                            \end{cases}$$
\begin{proof}
 We prove this by induction on $s := \rank_M(F)$. $s = 1$ is obvious, so assume $s > 1$. As $\rank(M \cdot H_G) = \rank(M)-1$, it is clear that $\rank_{M \cdot H_G}(F) \in \{\rank_M(F), \rank_M(F)-1\}.$ Hence we want to prove that $$\rank_{M \cdot H_G}(F) = \rank_M(F) -1 \textnormal{ if and only if }F \cup G = E.$$
 
 Assume $F \cup G = E$. Then by Lemma \ref{lemma_flats_of_product}, any flat $H$ covered by $F$ cannot be a flat of $M \cdot M_G$. Hence $\rank_{M \cdot H_G}(F) = \rank_M(F) -1$.
 
 Conversely, assume $\rank_{M \cdot H_G}(F) = \rank_M(F) -1$, but $F \cup G \subsetneq E$. Pick any flat $H$ of $M$ that is covered by $F$. If $H$ is a flat of $M \cdot H_G$, then we can use induction to see that $\rank_{M \cdot H_G}(F) = \rank_{M \cdot H_G}(H) + 1= \rank_M(H) + 1 = \rank_M(F)$, which is a contradiction. Hence $H$ is not a flat of $M \cdot H_G$. By Lemma \ref{lemma_flats_of_product}, there must be a flat $H'$ covering $H$ such that $H' \cup G = E$. In particular $\rank_M(H') = \rank_M(F)$ and flat axioms tell us that $F \cap H' = H$. Denote by $K = \textnormal{cl}_M(F \cup H')$. By semimodularity, we have $\rank_M(K) \leq \rank_M(F) + \rank_M(H) - \rank_M(F \cap H') = \rank_M(F)+1$. Thus, $K$ is a flat of $M$ covering $F$ such that $K \cup G = E$. But this implies that $F$ is not a flat of $M \cdot H_G$, which is a contradiction.
\end{proof}
\end{lemma}

\begin{prop}\label{prop_nested_intersection_zero}
 Let $M$ be a matroid on $[n]$ of rank at least 2 and $M_{\curly{G}} = H_{G_1} \cdot \dots \cdot H_{G_c}$ a nested matroid with $c < n - \corank(M)$. Then $M \cdot M_{\curly{G}} = 0$ if and only if the following holds: There exists an $i = 1,\dots,c$ and a flat $F$ of $M$ such that 
 \begin{enumerate}
  \item $\rank_M(F) = c - i + 1$.
  \item $F \cup G_i = E$.
 \end{enumerate}
 \begin{proof}
We prove this by induction on $c$. The case $c = 1$ is a reformulation of Lemma \ref{lemma_vanishing}. Now let $c > 1$. We write $M' := M \cdot H_{G_1}$ and $N' := H_{G_2} \cdot \dots \cdot H_{G_c}$. Hence $N'$ is a nested matroid of corank $c-1$.

First assume $M \cdot M_{\curly{G}} = 0$. Then either $M' = 0$ or $M' \cdot N' = 0$. If $M' = 0$, Lemma \ref{lemma_vanishing} tells us that there exists a rank one flat $F$ of $M$ such that $F \cup G_1 = E$ and we are done. If $M' \neq 0$, but $M' \cdot N' = 0$, we can use induction to see that there exists an $i' = 1,\dots,c-1$ and a flat $F$ of $M'$ such that
\begin{enumerate}
 \item $\rank_{M'}(F) = c-i'$.
 \item $F \cup G_{i'+1} = E$.
\end{enumerate}
As $M'$ is a quotient of $M$, $F$ is also a flat of $M$. We now distinguish two cases. If $\rank_M(F) = \rank_{M'}(F) + 1$, then by Lemma \ref{lemma_rank_down}, we must have $F \cup G_1 = E$ and since $\rank_M(F) = c-i'+1 \leq c$, we can replace $F$ by any flat of rank $c$ that contains it. If $\rank_M(F) = \rank_{M'}(F)$, we can pick $i = i'+1$ and the two conditions are fulfilled.

Conversely, assume we have a flat $F$ of $M$ and an $i \in \{1,\dots,c\}$ with the required properties. If $i = 1$, then $F$ is a flat of rank one in $M \cdot H_{G_1} \cdot \dots H_{G_{c-1}}$ by Lemmas \ref{lemma_flats_of_product} and \ref{lemma_rank_down}. So  $M \cdot M_{\curly{G}} = 0$ by Lemma \ref{lemma_vanishing}. Hence we assume $i > 1$. First assume that $F$ is not a flat of $M'$. If $F$ is not even a flat of $H_{G_1}$, then we must have $F \cup G_1 = E \wo \{x\}$ for some $x \in E$ by definition of $H_{G_1}$. The properties of the flats of a matroid dictate that there has to be a flat $F'$ of $M$  covering $F$ such that $x \in F'$, so $F' \cup G_1 = E$. If $F$ is a flat of $H_{G_1}$, then by Lemma \ref{lemma_flats_of_product}, there also exists such a flat $F'$. As $\rank_M(F') = \rank_M(F) + 1 = c-i+2 \leq c$, we can simply replace $F$ by a larger flat fulfilling the two conditions (1) and (2).

Thus we can assume that $F$ is indeed a flat of $M'$. Now we distinguish two cases: If $\rank_M(F) = \rank_{M'}(F) +1$, then $F \cup G_1 = E$ by Lemma \ref{lemma_rank_down} and since $G_1 \subseteq G_2$, we can apply induction to see that $M' \cdot N' = 0$. If $\rank_M(F) = \rank_{M'}(F)$, then we have $\rank_{M'}(F) = c -i +1 = (c-1) - (i-1) + 1$ and again we have $M' \cdot N' = 0$ by induction.
\end{proof}
\end{prop}

\subsection{Poincar\'e duality and the matroid polytope algebra}\label{subsection_poincare_duality}

By Theorem \ref{theorem_counting_nested} and the symmetry of the Eulerian numbers, we already know that the free $\Z$-modules $\M_{r,n}$ and $\M_{n-r+1,n}$ are isomorphic. From a purely matroid-theoretic perspective this might seem somewhat odd, as one would maybe rather have expected $\M_{r,n}$ to be isomorphic to $\M_{n-r,n}$. However, from a geometric perspective, the statement makes immediate sense, as the corresponding matroid fans have complementary dimensions. We will now see that there is another geometric interpretation of the ring $\M_n$, which makes it immediately clear what the isomorphism must be.

\begin{prop}\label{prop_cohomology}
 The ring $\M_n$ is isomorphic to the cohomology ring $A^*(X(\textnormal{Perm}_n))$ of the toric variety corresponding to the normal fan of the permutohedron of order $n$.
 \begin{proof}
  By \cite{fstoricintersection}, $A^*(X(\textnormal{Perm}_n))$ is the ring of all tropical cycles which are supported on some skeleton of the normal fan of the permutohedron. This normal fan is $B(U_{n,n})$, the fan of all chains, so we see that $\M_n$ is in fact a subring of $A^*(X(\textnormal{Perm}_n))$. It is a classical fact that the permutohedral variety has Eulerian numbers as Betti numbers, so the claim follows from Theorem \ref{theorem_counting_nested}.
 \end{proof}
\end{prop}

\begin{corollary}\label{theorem_poincare}
 Let $2 \leq n, 1 \leq r \leq n$. The intersection product on $\M_n$ induces a perfect pairing
 $$\M_{r,n} \times \M_{n-r+1,n} \to \M_{1,n} \cong \Z\enspace,$$
 i.e.\ it induces an isomorphism $\M_{r,n} \to \textnormal{Hom}(\M_{n-r+1,n},\Z) \cong \M_{n-r+1,n}.$
 \begin{proof}
  This is immediate from Proposition \ref{prop_cohomology}, as $X(\textnormal{Perm}_n)$ is a smooth and complete toric variety.
 \end{proof}
\end{corollary}

\begin{remark}
 A more general version of this can be found in \cite[Theorem 6.19]{ahkhodgetheory}, where Poincar\'e duality is shown for the ring of cycles supported on skeleta of an arbitrary Bergman fan.
\end{remark}

It has been shown \cite{fstoricintersection, jystableintersection} that the intersection ring $Z_{n-1}^{\textnormal{fan}}$ of all tropical fan cycles in $\tpn{n}$ is isomorphic to McMullen's polytope algebra $\Pi_{n-1}$ \cite{mpolytopealgebra} (considered over $\Q$). This is the algebra generated by symbols $[P]$ for each polytope in $\R^{n-1}$, modulo translations and the identity $[P \cup Q] = [P] + [Q] - [P \cap Q]$ whenever $P \cup Q$ is a polytope.

The isomorphism is defined by mapping the class $[P]$ of a polytope to
$$\exp(P) = \sum_{i=0}^{n-1} \frac{1}{i!} H_P^i\enspace,$$
where $H_P$ denotes the tropical hypersurface dual to $P$. As a set, this is just the codimension one skeleton of the normal fan of $P$. The weight of a maximal cell of $H_P$ is the lattice length of the dual edge of $P$. $H_P^i$ is the $i$-fold intersection product of the hypersurface.

It follows that $\M_n^\Q := \M_n \otimes_\Z \Q$ is isomorphic to a subring of $\Pi_{n-1}$. We can identify this subring precisely. For a matroid $M$ (possibly with loops), its \emph{matroid polytope} is 
$$P_M := \conv\left\{\sum_{i \in B} e_i;\; B \textnormal{ a basis of }M\right\} \subseteq \R^n\enspace.$$
Forgetting the last coordinate is a linear equivalence on $P_M$, so we can consider $[P_M]$ as an element of $\Pi_{n-1}$. Also, the normal fan of $P_M$ has a lineality space containing $\textbf{1}$, so its hypersurface can be seen as a cycle in $\tpn{n}$.

\begin{corollary}\label{corollary_polytope_algebra}
 Under the isomorphism $\exp$, we have
 $$\M_n^\Q \cong \Q[ [P_M]; M \textnormal{ a matroid on }\{1,\dots,n\}]\enspace.$$
Note that $M$ is allowed to have loops.
 \begin{proof}
  It is easy to see that for any matroid polytope $P_M$, its hypersurface is contained in the codimension one skeleton of the normal fan of the permutohedron. So by Proposition \ref{prop_cohomology} $\exp$ maps $[P_M]$ into $\M_n$. It remains to see that $\exp$ is surjective onto $\M_n^\Q$. By Corollary \ref{corollary_corank_one}, it suffices to show that every corank one matroid $H_G$ is in the image. But the matroid fan of $H_G$ is the hypersurface dual to the polytope $\conv\{e_i, i \notin G\}$, which is the matroid polytope of $(H_G)^*$. 
 \end{proof}
 \end{corollary}

 \begin{remark}\label{remark_derksen_fink}
  We already mentioned the matroid polytope modules studied by Derksen and Fink in \cite{dfvaluative}. The additive structure they consider is basically the same as on McMullen's polytope algebra, except that polytopes differing by a translation are not considered equal. In particular, it follows from Corollary \ref{corollary_polytope_algebra} that $\M_{r,n}$ is a quotient module of their module $P_M(n,r)$. One of their main results is the fact that $\curly{G}(M)$ is the universal valuative matroid invariant, in particular it is a linear map on $P_M(n,r)$. Hence Theorem \ref{theorem_g_invariant} would also follow from their result and the fact that $\curly{G}$ is compatible with translations of polytopes.
 \end{remark}

\section{Outlook}\label{section_outlook}

In this section we outline some interesting questions and connections for further research.

\subsection{Matroid homology}

For each $r \geq 1$ and each $i \in [n]$ there are two natural $\Z$-module homomorphisms $d_i: \M_{r,n} \to \M_{r,n-1}$ and $c_i: \M_{r,n} \to \M_{r-1,n-1}$ given on matroids by
\begin{align*}
 d_i(M) &:= \begin{cases}
             M\wo i, &\textnormal{if } i\textnormal{ is not a coloop of }M,\\
             0,&\textnormal{otherwise}.
            \end{cases}\\
 c_i(M) &:= \begin{cases}
             M/i, &\textnormal{if } \cl_M(\{i\}) = \{i\},\\
             0,&\textnormal{otherwise.}
            \end{cases}
\end{align*}

It is not obvious that these give well-defined maps. Tropical geometry comes to the rescue also in this case: Geometrically, $d_i$ corresponds to the \emph{push-forward} of tropical cycles along a coordinate projection and $c_i$ is the \enquote{intersection product with the hyperplane at $x_i = \infty$}. In both cases it is known that these operations commute with taking sums of cycles (see \cite{smatroidintersection} for proofs). These operations are defined to be zero whenever the dimension of the result does not match the expected dimension.

Let $\M = \oplus_{n \geq 0} \M_n$. In \cite{abgwhomology}, the authors show that one can define boundary maps on the free abelian group over all matroids using alternating sums of deletions or contractions. It is not hard to see that the same works for the maps $d_i$ and $c_i$, i.e.\ if we set
\begin{align*}
 \partial_d: &\M \to \M, M \mapsto \sum (-1)^i d_i(M)\\
 \partial_c: &\M \to \M, M \mapsto \sum (-1)^i c_i(M)\enspace,
\end{align*}
then $\partial_d^2 = \partial_c^2 = 0$. This allows us to define homology groups of matroids or minor-closed classed of matroids. Computational experiments suggest that when taking all matroids these homology groups always vanish, which seems not at all obvious.

\subsection{The matroid of matroids and the polytope of matroids}

Identifying each matroid with its indicator vector of chains $v_M$ makes the set of all loopfree matroids of rank $r$ into a matroid. So far, we have proven rather little about this matroid. We know its rank, which is $A_{r-1,n}$. Furthermore it is, by definition, realizable over any field of characteristic zero. 

Of particular interest are the circuits. More precisely, it would be interesting to understand the kernel of the map $\Phi_{r,n}: \M_{r,n}^\free \to V_{r,n}$, especially since the total number of loopfree matroids of rank $r$ on $n$ labeled elements is obviously $A_{r-1,n} + \dim \ker \Phi_{r,n}$.

If one considers the vector space $\M_n^\Q$ in the coordinates given by the basis of nested matroids, one can take the convex hull of the points corresponding to matroids. One can show that this is an empty lattice polytope (i.e.\ it has no interior lattice points), whose vertices are exactly the matroids. Since various matroid invariants define linear maps on $\M_n$, as shown in \ref{subsection_tutte}, this provides a new approach to open extremality questions, such as the one posed by Bonin and de Mier \cite{bmlatticeofcyclic} about the maximal possible number of cyclic flats of a matroid. Hence it seems essential to understand the combinatorial structure of this polytope of matroids.

\subsection{Regular subdivisions of matroid polytopes}

It was suggested in \cite{ahrrationalequivalence} that there should be an analogue of the polytope algebra which is isomorphic to the intersection ring of all tropical cycles (not just fans). The correct object would likely be an algebra of regular subdivisions of polytopes. One could again consider the subalgebra generated by matroid polytopes and their regular subdivisions and study its tropical counterpart. This might provide an interesting approach to understanding regular subdivisions of matroid polytopes better.

\bibliographystyle{halpha_with_href}
\bibliography{bibliography}

\end{document}